\documentclass{article} 
\usepackage{nips14submit_e,times}
\usepackage{hyperref}
\usepackage{array}
\usepackage{url}
\usepackage{booktabs}
\usepackage{tabu}
\usepackage{latexsym}
\usepackage{psfrag}
\usepackage{amsmath} 
\usepackage{amssymb}
\usepackage{bbm}
\usepackage{amsthm} 
\usepackage{mathrsfs} 
\usepackage{comment}
\usepackage{enumerate}
\usepackage{url}
\usepackage{graphicx}

\usepackage{rotating}
\usepackage[config, labelfont={normalsize,bf}, textfont=normalsize]{caption,subfig}
\usepackage{epsfig}
\usepackage[stable]{footmisc}
\usepackage{endnotes} 
\usepackage[english]{babel}
\usepackage{makecell}
\usepackage{rotating}
\usepackage{multirow}
\usepackage{cleveref}

\theoremstyle{definition}

\theoremstyle{plain}
\newtheorem{theorem}{Theorem}

\newtheorem{lemma}{Lemma}

\newtheorem{shortassumption}{A\!}

\theoremstyle{definition}

\theoremstyle{remark}

\newcommand{\R}{\mathbb{R}}

\renewcommand{\hat}[1]{\widehat{#1}}

\newcommand{\calI}{\mathcal{I}}

\newcommand{\tsum}{\textstyle\sum}

\renewcommand{\P}{\mathbb{P}}
\newcommand{\E}{\mathbb{E}}

\renewcommand{\L}{\mathcal{L}}

\newcommand{\ve}{\varepsilon}

\newcommand{\tr}{\operatorname{tr}}

\newcommand{\var}{\operatorname{var}}

\newcommand{\ttop}{^{\top}}

\newcommand{\op}{_{\text{op}}}

\newcommand{\hatbeta}{\hat{\beta}}

\newcommand{\hatbetarho}{\hat{\beta}_{\rho}}

\newcommand{\mspe}{\text{mspe}}

\newcommand\wcline[1]{%
    \noalign{\xdef\origarrayrulewidth{\the\arrayrulewidth}%
    \global\arrayrulewidth 4\arrayrulewidth}%
    \cline{#1}%
    \noalign{\global\arrayrulewidth\origarrayrulewidth}%
}

\newcommand\independent{\protect\mathpalette{\protect\independenT}{\perp}}
    \def\independenT#1#2{\mathrel{\setbox0\hbox{$#1#2$}%
    \copy0\kern-\wd0\mkern4mu\box0}} 

\newcommand{\1}{{\bf{1}}}
\newcommand{\de}{\overset{\mathcal{L}}{=}}
\newcommand{\ts}{\textstyle}
\newcommand{\mnorm}[1]{\left\vert\kern-1.5pt\left\vert\kern-1.5pt\left\vert #1\right\vert\kern-1.5pt\right\vert\kern-1.5pt\right\vert}

\title{A Residual Bootstrap for High-Dimensional Regression with Near Low-Rank Designs}

\author{
Miles E.~Lopes\\
Department of Statistics\\
University of California, Berkeley\\
Berkeley, CA 94720 \\
\texttt{mlopes@stat.berkeley.edu}
}

\nipsfinalcopy 

\begin{document}

\maketitle

\begin{abstract}
We study the residual bootstrap (RB) method in the context of high-dimensional linear regression. Specifically, we analyze the distributional approximation of linear contrasts $c\ttop (\hat{\beta}_{\rho}-\beta)$, where $\hat{\beta}_{\rho}$ is a ridge-regression estimator. When regression coefficients are estimated via least squares,  classical results show that RB consistently approximates the laws of contrasts, provided that $p\ll n$, where the design matrix is of size $n\times p$. Up to now, relatively little work has considered how additional structure in the linear model may extend the validity of RB to the setting where $p/n\asymp 1$. In this setting, we propose a version of RB that resamples residuals obtained from ridge regression. Our main structural assumption on the design matrix is that it is nearly low rank --- in the sense that its singular values decay according to a power-law profile. Under a few extra technical assumptions, we derive a simple criterion for ensuring that RB consistently approximates the law of a given contrast. We then specialize this result to study confidence intervals for mean response values $X_i\ttop \beta$, where $X_i\ttop$ is the $i$th row of the design.
More precisely, we show that conditionally on a Gaussian design with near low-rank structure,
RB \emph{simultaneously} approximates all of the laws $X_i\ttop(\hat{\beta}_{\rho}-\beta)$, $i=1,\dots,n$. This result is also notable as it imposes no sparsity assumptions on $\beta$. Furthermore, since our consistency results are formulated in terms of the Mallows (Kantorovich) metric, the existence of a limiting distribution is not required.
 \end{abstract}
\section{Introduction}\label{sec:intro}
Until recently, much of the emphasis in the theory of high-dimensional statistics has been on ``first order'' problems, such as estimation and prediction. As the understanding of these problems has become more complete, attention has begun to shift increasingly towards ``second order'' problems, dealing with hypothesis tests, confidence intervals, and uncertainty quantification~\cite{zhang2014confidence,javanmard2013hypothesis,javanmard2013confidence,buhlmannridge,vandegeerci2013,lee2013}. In this direction, much less is understood about the effects of structure, regularization, and  dimensionality --- leaving many questions open.
One collection of such questions that has attracted growing interest deals with the operating characteristics of the \emph{bootstrap} in high dimensions~\cite{lahiribootstrap2013,binyuridge2013,chernozhukovmultiplier}. Due to the fact that bootstrap is among the most widely used tools for approximating the sampling distributions of test statistics and estimators, there is much practical value in understanding what factors allow for the bootstrap to succeed in the high-dimensional regime.%

\paragraph{The regression model and linear contrasts.} In this paper, we focus our attention on high-dimensional linear regression, and our aim is to know when the residual bootstrap (RB) method consistently approximates the laws of \emph{linear contrasts}. (A review of  RB is given in Section~\ref{setupmethod}.)

To specify the model, suppose that we observe a response vector $Y\in \R^n$, generated according to
\begin{equation}\label{model}
Y=X \beta+\ve,
\end{equation}
where  $X\in\R^{n\times p}$ is the observed design matrix, $\beta \in\R^p$ is an unknown vector of coefficients, and the error variables $\ve=(\ve_1,\dots,\ve_n)$ are drawn i.i.d.~from an unknown distribution $F_0$, with mean $0$ and unknown variance $\sigma^2<\infty$. As is conventional in high-dimensional statistics, we assume the model~\eqref{model} is embedded in a sequence of models indexed by $n$. Hence, we allow $X$, $\beta$, and $p$ to vary implicitly with $n$. We will leave $p/n$ unconstrained until Section~\ref{ridgeconsist}, where we will assume $p/n\asymp 1$ in Theorem~\ref{ridgemspe}, and then in Section~\ref{sec:GaussianDesigns}, we will assume further that $p/n$ is bounded strictly between 0 and 1. The distribution $F_0$ is fixed with respect to $n$, and none of our results require $F_0$ to have more than four moments.

Although we are primarily interested in cases where the design matrix $X$ is deterministic, we will also study the performance of the bootstrap conditionally on a Gaussian design. For this reason, we will use the symbol $\E[\dots |X]$ even when the design is non-random so that confusion does not arise in relating different sections of the paper. Likewise, the symbol $\E[\dots]$ refers to unconditional expectation over all sources of randomness. Whenever the design is random, we will assume $X\independent \ve$, denoting the distribution of $X$ by $\P_X$, and the distribution of $\ve$ by $\P_{\ve}$.

Within the context of the regression, we will be focused on linear contrasts %
$c\ttop (\hat{\beta}-\beta)$, where $c\in \R^p$ is a fixed vector and $\hat{\beta}\in\R^p$ is an estimate of $\beta$. The importance of contrasts arises from the fact that they unify many questions about a linear model. For instance, testing the significance of the $i$th coefficient $\beta_i$ may be addressed by choosing $c$ to be the standard basis vector $c\ttop=e_i\ttop$.
Another important problem is quantifying the uncertainty of point predictions, which may be addressed by choosing  $c\ttop= X_i\ttop$, i.e.~the $i$th row of the design matrix. In this case, an approximation to the law of the contrast leads to a confidence interval for the mean response value $\E[Y_i]=X_i\ttop \beta$.
Further applications of contrasts occur in the broad topic of ANOVA~\cite{lehmann}.

\paragraph{Intuition for structure and regularization in RB.} The following two paragraphs explain the core conceptual aspects of the paper. To understand the role of regularization in applying RB to high-dimensional regression, it is helpful to think of RB in terms of two ideas. First, if $\hat{\beta}_{\text{LS}}$ denotes the ordinary least squares estimator, then it is a simple but important fact that contrasts can be written as $c\ttop (\hat{\beta}_{\text{LS}}-\beta) = a\ttop \ve$ where $a\ttop \!\!:=c\ttop (X\ttop X)^{-1}X\ttop$. Hence, if it were possible to sample directly from $F_0$, then the law of any such contrast could be easily determined. Since $F_0$ is unknown, the second key idea is to use the residuals of \emph{some} estimator $\hat{\beta}$ as a proxy for samples from $F_0$. When $p\ll n$, the least-squares residuals are a good proxy~\cite{freedman1981,bickel1983}. However, it is well-known that least-squares tends to overfit when $p/n\asymp 1$. When $\hat{\beta}_{\text{LS}}$ fits ``too well'', this means that its residuals are ``too small'', and hence they give a poor proxy for $F_0$. Therefore, by using a regularized estimator $\hat{\beta}$, overfitting can be avoided, and the residuals of $\hat{\beta}$ may offer a better way of obtaining ``approximate samples'' from $F_0$. 

The form of regularized regression we will focus on is \emph{ridge regression}:
\begin{equation}\label{ridgedef}
\hatbetarho:=(X\ttop X+\rho I_{p\times p})^{-1}X\ttop Y,
\end{equation}
where $\rho>0$ is a user-specificed regularization parameter.
 As will be seen in Sections~\ref{secmspelink} and~\ref{ridgeconsist}, the residuals obtained from ridge regression lead to a particularly good approximation of $F_0$ when the design matrix $X$ is nearly low-rank, in the sense that most of its singular values are close to 0. In essence, this condition is a form of sparsity, since it implies that the rows of $X$ nearly lie in a low-dimensional subspace of $\R^p$. However,  this type of structural condition has a significant advantage over the the more well-studied assumption that $\beta$ is sparse. Namely, the assumption that $X$ is nearly low-rank can be inspected directly in practice --- whereas sparsity in $\beta$ is typically unverifiable. In fact, our results will impose no conditions on $\beta$, other than that $\|\beta\|_2$ remains bounded as $(n,p)\to\infty$. Finally, it is worth noting that the occurrence of near low-rank design matrices is actually very common in applications,
and is often referred to as \emph{collinearity}~\cite[ch. 17]{draper1998}.

 \paragraph{Contributions and outline.} The primary contribution of this paper is a complement to the work of Bickel and Freedman~\cite{bickel1983} (hereafter B\&F 1983) --- who showed that in general, the RB method fails to approximate the laws of least-squares contrasts $c\ttop (\hat{\beta}_{\text{LS}}-\beta)$ when $p/n\asymp 1$. Instead, we develop an alternative set of results, proving that even when $p/n\asymp 1$, RB can successfully approximate the laws of ``ridged contrasts'' $c\ttop (\hat{\beta}_{\rho}-\beta)$ for many choices of $c\in \R^p$, provided that the design matrix $X$ is nearly low rank. A particularly interesting consequence of our work is that RB successfully approximates the law $c\ttop (\hat{\beta}_{\rho}-\beta)$ for a certain choice of $c$ that was shown in B\&F 1983 to ``break'' RB when applied to least-squares. Specifically, such a $c$ can be chosen as one of the rows of $X$ with a high \emph{leverage score} (see Section~\ref{sec:sims}). This example  corresponds to the practical problem of setting confidence intervals for mean response values $\E[Y_i]=X_i\ttop \beta$. (See~\cite[p. 41]{bickel1983}, as well as Lemma~\ref{propGaussian} and Theorem~\ref{lasthm} in Section~\ref{sec:GaussianDesigns}). Lastly, from a technical point of view, a third notable aspect of our results is that they are formulated in terms of the Mallows-$\ell_2$ metric, which  frees us from having to impose conditions that force a limiting distribution to exist. 
 
 Apart from B\&F 1983, the most closely related works we are aware of are the recent papers~\cite{lahiribootstrap2013} and~\cite{binyuridge2013}, which also consider RB in the high-dimensional setting. However, these works focus on role of sparsity in $\beta$ and do not make use of low-rank structure in the design, whereas our work deals only with structure in the design and imposes no sparsity assumptions on $\beta$.

 The remainder of the paper is organized as follows. In Section~\ref{setupmethod}, we formulate the problem of approximating the laws of contrasts, and describe our proposed methodology for RB based on ridge regression. Then, in Section~\ref{sec:main} we state several results that lay the groundwork for Theorem~\ref{lasthm}, which shows that that RB can successfully approximate all of the laws $\L(X_i\ttop(\hat{\beta}_{\rho}-\beta)|X)$, $i=1,\dots,n$, conditionally on a Gaussian design. Due to space constraints, all proofs are deferred to material that will appear in a separate work.

\paragraph{Notation and conventions.} If $U$ and $V$ are random variables, then $\L(U|V)$ denotes the law of $U$, conditionally on $V$. If $a_n$ and $b_n$ are two sequences of real numbers, then the notation $a_n\lesssim b_n$ means that there is an absolute constant $\kappa_0>0$ and an integer $n_0 \geq 1$ such that $a_n \leq \kappa_0 b_n$ for all $n\geq n_0$. The notation $a_n\asymp b_n$ means that $a_n \lesssim b_n$ and $b_n\lesssim a_n$. 
For a square matrix $A\in\R^{k\times k}$ whose eigenvalues are real, we denote
  them by $\lambda_{\min}(A)=\lambda_k(A)\leq \cdots \leq
  \lambda_1(A)=\lambda_{\max}(A)$.

\section{Problem setup and methodology}\label{setupmethod}

\paragraph{Problem setup.} For any $c\in\R^p$, it is clear that conditionally on $X$, the law of $c\ttop (\hatbetarho-\beta)$ is completely determined by $F_0$, and hence it makes sense to use the notation
\begin{equation}\label{psidef}
\Psi_{\rho}(F_0;c):= \L\big(c\ttop (\hatbetarho-\beta)\big\bracevert X\big).
\end{equation}
The problem we aim to solve is to approximate the distribution $\Psi_{\rho}(F_0;c)$ for suitable choices of $c$.

\paragraph{Review of the residual bootstrap (RB) procedure.} We briefly explain the steps involved in  the residual bootstrap procedure, applied to the ridge estimator $\hat{\beta}_{\rho}$ of $\beta$. To proceed somewhat indirectly, consider the following ``bias-variance'' decomposition of $\Psi_{\rho}(F_0;c)$, conditionally on $X$,
\begin{equation}\label{decompose}
\Psi_{\rho}(F_0;c) =  \underbrace{\L\big(c\ttop \big(\hatbetarho-\E[\hatbetarho|X]\big)\big\bracevert X\big)}_{=: \ \Phi_{\rho}(F_0;c)} \ + \ \underbrace{c\ttop\big(\E[\hatbetarho|X]-\beta\big)}_{=: \, \text{bias}(\Phi_{\rho}(F_0;c))}.
\end{equation}
Note that the distribution $\Phi(F_0;c)$ has mean zero, and so that the second term on the right side is the bias of $\Phi_{\rho}(F_0;c)$ as an estimator of $\Psi_{\rho}(F_0;c)$. Furthermore, the distribution $\Phi_{\rho}(F_0;c)$ may be viewed as the ``variance component'' of $\Psi_{\rho}(F_0;c)$.
We will be interested in situations where the regularization parameter $\rho$ may be chosen small enough so that the bias component is small. In this case, one has
$ \Psi_{\rho}(F_0;c) \approx \Phi_{\rho}(F_0;c),$
and then it is enough to find an approximation to the law $\Phi_{\rho}(F_0;c)$, which is unknown. To this end, a simple manipulation of $c\ttop(\hat{\beta}_{\rho}-\E[\hat{\beta}_{\rho}])$ leads to
\begin{equation}\label{PhiF0}
\Phi_{\rho}(F_0;c)=\mathcal{L}(c\ttop (X\ttop X+\rho I_{p\times p})^{-1}X\ttop \ve\big\bracevert X).
\end{equation}

Now, to approximate $\Phi_{\rho}(F_0;c)$, let $\hat{F}$ be any centered estimate of $F_0$. (Typically, $\hat{F}$ is obtained by using the centered residuals of some estimator of $\beta$, but this is not necessary in general.)
Also, let $\ve^*= (\ve_1^*,\dots,\ve_n^*)\in \R^n$ be an i.i.d. sample from $\hat{F}$. Then, replacing $\ve$ with $\ve^*$ in line~\eqref{PhiF0} yields
\begin{equation}\label{phihat}
\Phi_{\rho}(\hat{F};c)=\mathcal{L}(c\ttop(X\ttop X+\rho I_{p\times p})^{-1}X\ttop \ve^*\big\bracevert X).
\end{equation}
At this point, we define the (random) measure $\Phi_{\rho}(\hat{F};c)$ to be the RB approximation to $\Phi_{\rho}(F_0;c)$. Hence, it is clear that the RB approximation is simply a ``plug-in rule''.

\paragraph{A two-stage approach.} An important feature of the procedure just described is that we are free to use any centered estimator $\hat{F}$ of $F_0$.
This fact offers substantial flexibility in approximating $\Psi_{\rho}(F_0;c)$. One way of exploiting this flexibility is to consider a two-stage approach, where a ``pilot'' ridge estimator $\hat{\beta}_{\varrho}$ is used to first compute residuals whose centered empirical distribution function is $\hat{F}_{\varrho}$, say.
 Then, in the second stage, the distribution $\hat{F}_{\varrho}$ is used to approximate $\Phi_{\rho}(F_0;c)$ via the relation~\eqref{phihat}.
To be more detailed, if $(\hat{e}_1(\varrho),\dots,\hat{e}_n(\varrho))=\hat{e}(\varrho) :=Y-X\hat{\beta}_{\varrho}$ are the residuals of $\hat{\beta}_{\varrho}$, then we define $\hat{F}_{\varrho}$ to be the distribution that places mass $1/n$ at each of the values $\hat{e}_i(\varrho)-\bar{e}(\varrho)$ with $\bar{e}(\varrho):=\frac{1}{n}\sum_{i=1}^n \hat{e}_i(\varrho)$. Here, it is important to note that the value $\varrho$ is chosen to optimize $\hat{F}_{\varrho}$ as an approximation to $F_0$. By contrast, the choice of $\rho$ depends on the relative importance of width and coverage probability for confidence intervals based on $\Phi_{\rho}(\hat{F}_{\varrho};c)$. Theorems~\ref{basicprop},~\ref{ridgemspe}, and~\ref{lasthm} will offer some guidance in selecting $\varrho$ and $\rho$.

\paragraph{Resampling algorithm.} 
 To summarize the discussion above, if $B$ is user-specified number of bootstrap replicates, our proposed method for approximating $\Psi_{\rho}(F_0;c)$ is given below.
\begin{enumerate}\itemsep 0.07cm
\item Select $\rho$ and $\varrho$, and compute the residuals $\hat{e}(\varrho)=Y-X\hat{\beta}_{\varrho}$.
\item Compute the centered distribution function $\hat{F}_{\varrho}$, putting mass $1/n$ at each $\hat{e}_i(\varrho)-\bar{e}(\varrho)$.
\item For $j=1,\dots, B$:
\begin{itemize}
\item Draw a vector $\ve^*\in \R^n$ of $n$ i.i.d. samples from $\hat{F}_{\varrho}$.
\item Compute $z_j:=c\ttop (X\ttop X+\rho I_{p\times p})^{-1}X\ttop \ve^*$.
\end{itemize}
\item Return the empirical distribution of $z_1,\dots,z_B$.
\end{enumerate}
Clearly, as $B\to\infty$, the empirical distribution of $z_1,\dots,z_B$ converges weakly to $\Phi_{\rho}(\hat{F}_{\varrho};c)$, with probability 1. 
As is conventional, our theoretical analysis in the next section will ignore Monte Carlo issues, and address only the performance of $\Phi_{\rho}(\hat{F}_{\varrho};c)$ as an approximation to $\Psi_{\rho}(F_0;c)$. 

\section{Main results}\label{sec:main}
The following metric will be central to our theoretical results, and has been a standard tool in the analysis of the bootstrap, beginning with the work of Bickel and Freedman~\cite{bickel1981}.
\paragraph{The Mallows (Kantorovich) metric.} 
For two random vectors $U$ and $V$ in a Euclidean space, the Mallows-$\ell_2$ metric is defined by
\begin{equation}
d_2^2(\L(U),\L(V)) := \inf_{\pi\in \Pi} \Big\{ \E\Big[\|U-V\|_2^2\Big] : (U,V)\sim \pi \Big\}
\end{equation}
where the infimum is over the class $\Pi$ of joint distributions $\pi$ whose marginals are $\L(U)$ and $\L(V)$. It is worth noting that convergence in $d_2$ is strictly stronger than weak convergence, since it also requires convergence of second moments. Additional details may be found in the paper~\cite{bickel1981}.

\subsection{A bias-variance decomposition for bootstrap approximation}
To give some notation for analyzing the bias-variance decomposition of $\Psi_{\rho}(F_0;c)$ in line~\eqref{decompose}, we define the following quantities based upon the ridge estimator $\hat{\beta}_{\rho}$. Namely, the variance is
$$
v_{\rho}=
v_{\rho}(X;c) := \var(\Psi_{\rho}(F_0;c)|X)
= \sigma^2 \|c\ttop (X\ttop X+\rho I_{p\times p})^{-1}X\ttop\|_2^2.
$$
To express the bias of $\Phi_{\rho}(F_0;c)$, we define the vector $\delta(X)\in\R^p$ according to
\begin{equation}
\textstyle 
\delta(X):= \beta-\E[\hatbetarho] = \big[I_{p\times p}-(X\ttop X+\rho I_{p\times p})^{-1}X\ttop X \big]\beta,
\end{equation}
and then put
\begin{align}
\textstyle
b^2_{\rho}=b^2_{\rho}(X;c):=\text{bias}^2(\Phi_{\rho}(F_0;c))
= (c\ttop \delta(X))^2.
\end{align}
We will sometimes omit the arguments of $v_{\rho}$ and $b_{\rho}^2$ to lighten notation. Note that $v_{\rho}(X;c)$ does not depend on $\beta$, and $b_{\rho}^2(X;c)$ only depends on $\beta$ through $\delta(X)$.

 The following result gives a regularized and high-dimensional extension of  some lemmas in Freedman's early work~\cite{freedman1981} on RB for least squares. The result does not require any structural conditions on the design matrix, or on the true parameter $\beta$.
 \begin{theorem}[consistency criterion]\label{basicprop}
 Suppose $X\in \R^{n\times p}$ is fixed. Let $\hat{F}$ be any estimator of $F_0$, and let $c\in \R^p$ be any vector such that $v_{\rho}=v_{\rho}(X;c)\neq 0$. Then with $\P_{\ve}$-probability 1, the following inequality holds for every $n\geq 1$, and every $\rho>0$,
 \begin{equation}\label{basicbound}
 d_2^2\Big(\ts\frac{1}{\sqrt{v_{\rho}}}\Psi_{\rho}(F_0;c),\ts\frac{1}{\sqrt{v_{\rho}}}\Phi_{\rho}(\hat{F};c)\Big) \leq 
\ts\frac{1}{\sigma^2} d_2^2(F_0,\hat{F})+\frac{b_{\rho}^2}{v_{\rho}}. \end{equation}
 \end{theorem}
 \paragraph{Remarks.} Observe that the normalization $1/\sqrt{v_{\rho}}$ ensures that the bound is non-trivial, since the distribution $\Psi_{\rho}(F_0;c)/\sqrt{v_{\rho}}$ has variance equal to 1 for all $n$ (and hence does not become degenerate for large $n$). To consider the choice of $\rho$, it is simple to verify that the ratio $b_{\rho}^2/v_{\rho}$ decreases monotonically as $\rho$ decreases. Note also that as $\rho$ becomes small, the variance $v_{\rho}$ becomes large, and likewise, confidence intervals based on $\Phi_{\rho}(\hat{F};c)$ become wider. In other words, there is a trade-off between the width of the confidence interval and the size of the bound~\eqref{basicbound}.

 \paragraph{Sufficient conditions for consistency of RB.} An important practical aspect of Theorem~\ref{basicprop} is that for any given contrast $c$, the variance $v_{\rho}(X;c)$ can be easily estimated, since it only requires an estimate of $\sigma^2$, which can be obtained from $\hat{F}$. Consequently, whenever theoretical bounds on $d_2^2(F_0,\hat{F})$ and $b_{\rho}^2(X;c)$ are available, the right side of line~\eqref{basicbound} can be controlled. In this way, Theorem~\ref{basicprop} offers a simple route for guaranteeing that RB is consistent. In Sections~\ref{secmspelink} and~\ref{ridgeconsist} to follow, we derive a bound on $\E[d_2^2(F_0,\hat{F})|X]$ in the case where $\hat{F}$ is chosen to be $\hat{F}_{\varrho}$.
 Later on in Section~\ref{sec:GaussianDesigns}, we study RB consistency in the context of prediction with a Gaussian design, and there we derive high probability bounds on both $v_{\rho}(X;c)$ and $b_{\rho}^2(X;c)$ where $c$ is a particular row of $X$.

\subsection{A link between bootstrap consistency and MSPE}\label{secmspelink}
If $\hat{\beta}$ is an estimator of $\beta$, its mean-squared prediction error (MSPE), conditionally on $X$, is defined as
\begin{equation}\label{mspedef}
\text{mspe}(\hat{\beta}\,|X) :=\ts\frac{1}{n}\E\big[ \|X(\hat{\beta}-\beta)\|_2^2 \big\bracevert X\big].
\end{equation}
The previous subsection showed that in-law approximation of contrasts is closely tied to the approximation of $F_0$. We now take a second step of showing that if the centered residuals of an estimator $\hat{\beta}$ are used to approximate $F_0$, then the quality of this approximation can be bounded naturally in terms of  $\mspe(\hat{\beta}\,|X)$. 
This result applies to any estimator $\hat{\beta}$ computed from the observations~\eqref{model}.

\begin{theorem}\label{mspelink}
 Suppose $X\in\R^{n\times p}$ is fixed. Let $\hat{\beta}$ be any estimator of $\beta$, and let $\hat{F}$ be the empirical distribution of the centered residuals of $\hat{\beta}$. Also, let $F_n$ denote the empirical distribution of $n$ i.i.d. samples from $F_0$. Then for every $n\geq 1$,
\begin{equation}\label{predcons}
\E \,\big[ d_2^2(\hat{F},F_0)\big\bracevert X\big] \leq 2\,\text{\emph{mspe}}(\hat{\beta}\,|X)+2\,\E[d_2^2(F_n,F_0)]+\ts\frac{2\sigma^2}{n}.
\end{equation}
\end{theorem}

\paragraph{Remarks.} As we will see in the next section, the MSPE of ridge regression can be bounded in a sharp way when the design matrix is approximately low rank, and there we will analyze $\mspe(\hat{\beta}_{\varrho}|X)$ for the pilot estimator. Consequently, when near low-rank structure is available, the only remaining issue in controlling the right side of line~\eqref{predcons} is to bound the quantity $\E[d_2^2(F_n,F_0)|X] $. The very recent work of Bobkov and Ledoux~\cite{bobkovledoux} provides an in-depth study of this question, and they derive a variety bounds under different tail conditions on $F_0$.  We summarize one of their results below.
\begin{lemma}[Bobkov and Ledoux, 2014]\label{bobkovledouxmoment} If $F_0$ has a finite fourth moment, then
\begin{equation}
\E[d_2^2(F_n,F_0)]\lesssim \log(n)n^{-1/2}.
\end{equation}

\end{lemma}

%
%
%
%
{\bf{Remarks.}} The fact that the \emph{squared} distance is bounded at the rate of $\log(n)n^{-1/2}$ is an indication that $d_2$ is a rather strong metric on distributions.  For a detailed discussion of this result, see Corollaries 7.17 and 7.18 in the paper~\cite{bobkovledoux}.  Although it is possible to obtain faster rates when more stringent tail conditions are placed on $F_0$, we will only need a fourth moment, since the $\mspe(\hat{\beta}|X)$ term in Theorem~\ref{mspelink} will often have a slower rate than $\log(n) n^{-1/2}$, as discussed in the next section.

\subsection{Consistency of ridge regression in MSPE for near low rank designs}\label{ridgeconsist}
In this subsection, we show that when the tuning parameter $\varrho$ is set at a suitable rate, the pilot ridge estimator $\hat{\beta}_{\varrho}$ is consistent in MSPE when the design matrix is near low-rank ---  even when $p/n$ is large, and without any sparsity constraints on $\beta$. We now state some assumptions. 

\begin{shortassumption}\label{Ascaleeigs} 
There is a number $\nu>0$, and absolute constants $\kappa_1,\kappa_2>0$, such that
$$\kappa_1i^{-\nu} \leq \lambda_i(\hat{\Sigma}) \leq \kappa_2 i^{-\nu} \text{  \ \ \ \ \ \ \ for all \ \ \ } i=1,\dots,n\wedge p.$$
\end{shortassumption}

\begin{shortassumption}\label{Ascalerhofirst}
There are absolute constants $\theta,\gamma >0$, such that for every $n\geq 1$,
$\frac{\varrho}{n}=n^{-\theta}$ and $\frac{\rho}{n}=n^{-\gamma}$.
\end{shortassumption}
\begin{shortassumption}\label{Abetamag}
The vector $\beta\in \R^p$ satisfies $\|\beta\|_2\lesssim 1$.
\end{shortassumption}

 Due to Theorem~\ref{mspelink}, the following bound shows that the residuals of $\hat{\beta}_{\varrho}$ may be used to extract a consistent approximation to  $F_0$. Two other notable features of the bound are that it is \emph{non-asymptotic} and \emph{dimension-free}.

\begin{theorem}\label{ridgemspe}
Suppose that $X\in\R^{n\times p}$ is fixed and that Assumptions~\ref{Ascaleeigs}--\ref{Abetamag} hold, with $p/n\asymp 1$. Assume further that $\theta$ is chosen as $\theta=\frac{2\nu}{3}$ when $\nu\in (0,\ts\frac{1}{2})$, and  $\theta=\frac{\nu}{\nu+1}$ when $\nu>\ts\frac{1}{2}$. Then, 
\begin{equation}\label{mspebound2}
\text{\emph{mspe}}(\hat{\beta}_{\varrho}|X) \lesssim
\begin{cases}
& \,n^{-\frac{2\nu}{3}} \text{ \ \ \ \ \ if } \  \ \nu \in (0,\ts\frac{1}{2}),\\
&  \, n^{-\frac{\nu}{\nu+1}} \, \text{ \ \  if } \ \ \nu>\ts\frac{1}{2}.
\end{cases}
\end{equation}
Also, both bounds in~\eqref{mspebound2} are tight  in the sense that $\beta$ can be chosen so that $\hat{\beta}_{\varrho}$ attains either rate.
\end{theorem}
\paragraph{Remarks.} 
Since the eigenvalues $\lambda_i(\hat{\Sigma})$ are observable, they may be used to estimate $\nu$ and guide the selection of $\varrho/n=n^{-\theta}$.
However, from a practical point of view, we found it easier to select $\varrho$ via cross-validation in numerical experiments, rather than via an estimate of $\nu$.
\paragraph{A link with Pinsker's Theorem.} In the particular case when $F_0$ is a centered Gaussian distribution, the ``prediction problem'' of estimating $X\beta$ is very similar to estimating the mean parameters of a Gaussian sequence model, with error measured in the $\ell_2$ norm.  In the alternative sequence-model format, the decay condition on the eigenvalues of $\ts\frac{1}{n}X\ttop X$ translates into an ellipsoid constraint on the mean parameter sequence~\cite{tsybakov2009,wassermanallnp}. For this reason, Theorem~\ref{ridgemspe} may be viewed as ``regression version'' of $\ell_2$ error bounds for the sequence model under an ellipsoid constraint (cf. Pinsker's Theorem,~\cite{tsybakov2009,wassermanallnp}). Due to the fact that the latter problem has a very well developed literature, there may be various ``neighboring results'' elsewhere. Nevertheless, we could not find a direct reference for our stated MSPE bound in the current setup. For the purposes of our work in this paper, the more important point to take away from Theorem~\ref{ridgemspe} is that it can be coupled with Theorem~\ref{mspelink} for proving consistency of RB.

\subsection{Confidence intervals for mean responses, conditionally on a Gaussian design}\label{sec:GaussianDesigns}
In this section, we consider the situation where the design matrix $X$ has rows $X_i\ttop\in \R^p$ drawn i.i.d.~from a multivariate normal distribution $N(0,\Sigma)$, with $X\independent \ve$. (The covariance matrix $\Sigma$ may vary with $n$.)
Conditionally on a realization of $X$, we analyze the RB approximation of the laws $\Psi_{\rho}(F_0;X_i)=\L(X_i\ttop(\hatbetarho-\beta)|X)$. As discussed in Section~\ref{sec:intro}, this corresponds to the problem of setting confidence intervals for the mean responses $\E[Y_i]=X_i\ttop \beta$. 
Assuming that the population eigenvalues $\lambda_i(\Sigma)$ obey a decay condition, we show below in Theorem~\ref{lasthm} that RB succeeds with high $\P_X$-probability. Moreover, this consistency statement holds for all of the laws
 $\Psi_{\rho}(F_0;X_i)$ \emph{simultaneously}. That is, among the $n$ distinct laws $\Psi_{\rho}(F_0;X_i)$, $i=1,\dots,n$, even the worst bootstrap approximation is still consistent. 
We now state some population-level assumptions.
\begin{shortassumption}\label{Aop}
The operator norm of $\Sigma\in\R^{p\times p}$ satisfies $\|\Sigma\|_{\text{\emph{op}}} \lesssim 1$.
\end{shortassumption}
Next, we impose a decay condition on the eigenvalues of $\Sigma$. This condition also ensures that $\Sigma$ is invertible for each fixed $p$ --- even though the bottom eigenvalue may become arbitrarily small as $p$ becomes large. It is important to notice that we now use $\eta$ for the decay exponent of the population eigenvalues, whereas we used $\nu$ when describing the sample eigenvalues in the previous section.
\begin{shortassumption}\label{AeigSigma}
There is a number $\eta>0$, and absolute constants $k_1,k_2>0$, such that  for all $i=1,\dots,p$,
$$k_1i^{-\eta} \leq \lambda_i(\Sigma) \leq k_2 i^{-\eta}.$$
\end{shortassumption}
\begin{shortassumption}\label{Adim}
There are absolute constants $k_3,k_4\in(0,1)$ such that for all $n\geq 3$, we have the bounds\\ \indent \ \ \ \ \ \ \  \ $k_3\leq \ts\frac{p}{n}\leq k_4$ and  $p\leq n-2$.
\end{shortassumption}

The following lemma collects most of the effort needed in proving our final result in Theorem~\ref{lasthm}. Here it is also helpful to recall the notation $\rho/n= n^{-\gamma}$ and $\varrho/n=n^{-\theta}$ from Assumption~\ref{Ascalerhofirst}.
\begin{lemma}\label{propGaussian}
Suppose that the matrix $X\in \R^{n\times p}$ has rows $X_i\ttop$ drawn i.i.d. from $N(0,\Sigma)$, and that
Assumptions~\ref{Ascalerhofirst}--\ref{Adim}
hold. Furthermore, assume that $\gamma$ chosen so that $0<\gamma< \min\{\eta,1\}$.  Then, the statements below are true.

\noindent (i) (bias inequality) \\
Fix any $\tau>0$. Then, there is an absolute constant $\kappa_0>0$, such that for all large $n$, the following event holds with $\P_X$-probability at least \mbox{$1-n^{-\tau}-ne^{-n/16}$,}
\begin{equation}\label{biasineq}
\max_{1\leq i\leq n} b_{\rho}^2(X;X_i)\ \leq \kappa_0 \cdot n^{-\gamma} \cdot \ts (\tau+1)\log(n+2).
\end{equation}
~\\[-0.15cm]
\noindent (ii) (variance inequality)\\
There are absolute constants $\kappa_1,\kappa_2>0$ such that for all large $n$, the following event holds with $\P_X$-probability at least $1- 4n\exp(-\kappa_1 n^{\frac{\gamma}{\eta}})$,
\begin{equation}\label{varianceineq}
\max_{1\leq i\leq n}\ts \frac{1}{v_{\rho}(X;X_i)}  \leq  \kappa_2  n^{1-\frac{\gamma}{\eta}}.
\end{equation}
~\\[-0.2cm]
\noindent (iii) (mspe inequalities)\\
Suppose that $\theta$ is chosen as $\theta=2\eta/3$ when $\eta\in(0,\ts\frac{1}{2})$, and that $\theta$ is chosen as $\theta=\ts\frac{\eta}{1+\eta}$ when $\eta>\ts\frac{1}{2}$. Then, there are absolute constants $\kappa_3,\kappa_4,\kappa_5,\kappa_6>0$ such that for all large $n$,
\begin{equation*}\label{mspecombine}
\text{\emph{mspe}}(\hat{\beta}_{\varrho}|X)\leq 
\begin{cases}
& \kappa_4 n^{-\frac{2\eta}{3}} \text{ \ \ \ \ \ with $\P_X$-probability at least } \ \ 1-\exp(-\kappa_3 n^{2-4\eta/3}), \ \text{ if }  \ \eta \in (0,\ts\frac{1}{2})\\
& \kappa_6 n^{-\frac{\eta}{\eta+1}} \text{ \ \ \ with $\P_X$-probability at least } \ \ 1-\exp(-\kappa_5n^{\frac{2}{1+\eta}}), \ \ \ \ \ \text{ if } \ \  \eta>\ts\frac{1}{2}.
\end{cases}
\end{equation*}
\end{lemma}


\vspace{-0.3cm}
\paragraph{Remarks.} 
Note that the two rates in part (iii) coincide as $\eta$ approaches $1/2$. At a conceptual level, the entire lemma may be explained in relatively simple terms. 
Viewing the quantities $\mspe(\hat{\beta}_{\varrho}|X)$, $b_{\rho}^2(X;X_i )$  and $v_{\rho}(X;X_i)$ as functionals of a Gaussian matrix, the proof involves deriving concentration bounds for each of them. Indeed, this is plausible given that these quantities are smooth functionals of $X$. However, the difficulty of the proof arises from the fact that they are also highly non-linear functionals of $X$. 
\vspace{-0.05cm}

We now combine Lemmas~\ref{bobkovledouxmoment} and~\ref{propGaussian} with Theorems~\ref{basicprop} and~\ref{mspelink} to show that all of the laws $\Psi_{\rho}(F_0;X_i)$ can be simultaneously approximated via our two-stage RB method.

\begin{theorem}\label{lasthm}
Suppose that $F_0$ has a finite fourth moment, Assumptions \ref{Ascalerhofirst}--\ref{Adim} hold, and $\gamma$ is chosen so that $\frac{\eta}{1+\eta}< \gamma < \min\{\eta,1\}$. Also suppose that $\theta$ is chosen as $\theta=2\eta/3$ when
 $\eta\in(0,\ts\frac{1}{2})$, and  $\theta=\frac{\eta}{\eta+1}$ when  $\eta>\frac{1}{2}$.
 Then, there is a sequence of positive numbers $\delta_n$ with $\lim_{n\to\infty}\delta_n =0$, such that the event
\begin{equation}\label{lasthm1}
\E\Big[ \max_{1\leq i\leq n} d_2^2\Big(\ts\frac{1}{\sqrt{v_{\rho}}}\Psi_{\rho}(F_0; X_i), \,\ts\frac{1}{\sqrt{v_{\rho}}} \Phi_{\rho}(\hat{F}_{\varrho};X_i)\Big)\Big\bracevert X\Big] 
\leq \delta_n
\end{equation}
has $\P_X$-probability tending to 1 as $n\to\infty$.
\end{theorem}
\paragraph{Remark.} Lemma~\ref{propGaussian} gives explicit bounds on the numbers $\delta_n$, as well as the probabilities of the corresponding events, but we have stated the result in this way for the sake of readability.
%
%
%
\section{Simulations}\label{sec:sims}
In four different settings of $n,p,$ and the decay parameter $\eta$, we compared the nominal $90\%$ confidence intervals (CIs) of four methods: ``oracle'', ``ridge'', ``normal'', and ``OLS'', to be described below. In each setting, we generated $N_1:=100$ random designs $X$ with i.i.d. rows drawn from $N(0,\Sigma)$, where $\lambda_j(\Sigma)=j^{-\eta}$, $j=1,\dots,p$, and the eigenvectors of $\Sigma$ were drawn randomly by setting them to be the $Q$ factor in a $QR$ decomposition of a standard $p\times p$ Gaussian matrix. Then, for each realization of $X$, we generated $N_2:=1000$ realizations of $Y$ according to the model~\eqref{model}, where $\beta = \1/\|\1\|_2 \in \R^p$, and $F_0$ is the centered $t$ distribution on $5$ degrees of freedom, rescaled to have standard deviation $\sigma=0.1$. For each $X$, and each corresponding $Y$, we considered the problem of setting a $90\%$ CI for the mean response value $X_{i^{\star}}\ttop \beta$, where $X_{i^{\star}}\ttop$ is the row with the highest leverage score, i.e. $i^{\star} = \text{argmax}_{1\leq i\leq n} \,H_{ii}$ and $H:=X(X\ttop X)^{-1}X\ttop$. This problem was shown in B\&F 1983 to be a case where the standard RB method based on least-squares fails when $p/n\asymp 1$. Below, we refer to this method as ``OLS''.  

To describe the other three methods, ``ridge'' refers to the interval  $[X_{i^{\star}}\ttop \hat{\beta}_{\rho}-\hat{q}_{0.95},X_{i^{\star}}\ttop \hat{\beta}_{\rho}-\hat{q}_{0.05}]$, where $\hat{q}_{\alpha}$ is the $\alpha \%$ quantile of the numbers $z_1,\dots,z_B$ computed in the proposed algorithm in Section~\ref{setupmethod}, with $B=1000$ and $c\ttop =X_{i^{\star}}\ttop$. To choose the parameters $\rho$ and $\varrho$ for a given $X$ and $Y$, we first computed $\hat{r}$ as the value that optimized the MSPE error of a ridge estimator $\hat{\beta}_{r}$ with respect to 5-fold cross validation; i.e.~cross validation was performed for every distinct pair $(X,Y)$.
We then put $\varrho =5\hat{r}$ and $\rho= 0.1\hat{r}$, as we found the prefactors 5 and 0.1 to work adequately across various settings. (Optimizing $\varrho$ with respect to MSPE is motivated by Theorems~\ref{basicprop},~\ref{mspelink}, and~\ref{ridgemspe}. Also, choosing $\rho$ to be somewhat smaller than $\varrho$ conforms with the constraints on $\theta$ and $\gamma$ in Theorem~\ref{lasthm}.)
The method ``normal'' refers to the CI based on the normal approximation $\L(X_{i^{\star}}\ttop (\hat{\beta}_{\rho}-\beta)|X)\approx N(0,\hat{\tau}^2)$, where $\hat{\tau}^2 = \hat{\sigma}^2\|X_{i^{\star}}\ttop(X\ttop X+\rho I_{p\times p})^{-1}X\ttop\|_2^2$, $\rho=0.1\hat{r}$, and $\hat{\sigma}^2$ is the usual unbiased estimate of $\sigma^2$ based on OLS residuals.
The ``oracle'' method refers to the interval $[X_{i^{\star}}\ttop \hat{\beta}_{\rho}-\tilde{q}_{0.95},X_{i^{\star}}\ttop \hat{\beta}_{\rho}-\tilde{q}_{0.05}]$, with $\rho=0.1\hat{r}$, and $\tilde{q}_{\alpha}$ being the empirical $\alpha\%$ quantile of $X_i\ttop (\hat{\beta}_{\rho}-\beta)$ over all $1000$ realizations of $Y$ based on a given $X$. (This accounts for the randomness in $\rho=0.1\hat{r}$.)

Within a given setting of the triplet $(n,p,\eta)$, we refer to the ``coverage'' of a method as the fraction of the $N_1\times N_2=10^5$ instances where the method's CI contained the parameter $X_{i^{\star}}\ttop\beta$. Also, we refer to ``width'' as the average width of a method's intervals over all of the $10^5$ instances. The four settings of $(n,p,\eta)$ correspond to moderate/high dimension and moderate/fast decay of the eigenvalues $\lambda_i(\Sigma)$. Even in the moderate case of $p/n=0.45$,  the results show that the OLS intervals are too narrow and have coverage noticeably less than 90\%. As expected, this effect becomes more pronounced when $p/n=0.95$. The ridge and normal intervals perform reasonably well across settings, with both performing much better than OLS.  However, it should be emphasized that our study of RB is motivated by the desire to gain insight into the behavior of the bootstrap in high dimensions --- rather than trying to outperform particular methods. In future work, we plan to investigate the relative merits of the ridge and normal intervals in greater detail.

\begin{table}[h!]
\caption{Comparison of nominal $90\%$ confidence intervals} 
\centering 
\setlength\extrarowheight{3pt}
\begin{tabular}{cc|c|c|c|c|l}
\cline{3-6}
\cline{3-6}
& & \ \ \ \  oracle \ \ \ \  & \ \ \ ridge \ \ \ & \ \ \ normal \ \ \ & \ \ \ OLS \ \ \ \\\cline{1-6}
\multicolumn{1}{ |c  }{setting 1}  &
\multicolumn{1}{ |c| }{ \ \ \ \  width \   \ \ \ } & 0.21 & 0.20 & 0.23 & 0.16 & \\\cline{2-6}
\multicolumn{1}{ |c  }{$n=100, \ p=45, \ \ \eta=0.5$}                        &
\multicolumn{1}{ |c| }{coverage} & 0.90 & 0.87 & 0.91 & 0.81 &  \\ \wcline{1-6}
\multicolumn{1}{ |c  }{setting 2}  &
\multicolumn{1}{ |c| }{width} &  0.22 & 0.26 & 0.26 & 0.06 &   \\ \cline{2-6}
\multicolumn{1}{ |c  }{$n=100, \ p=95, \ \ \eta=0.5$}                        &
\multicolumn{1}{ |c| }{coverage} &  0.90 & 0.88 & 0.88 & 0.42 &    \\ 
\wcline{1-6}
\multicolumn{1}{ |c  }{setting 3}  &
\multicolumn{1}{ |c| }{width} &  0.20 & 0.21 & 0.22 & 0.16 &    \\ \cline{2-6}
\multicolumn{1}{ |c  }{$n=100, \ p=45, \ \ \eta=1$}                        &
\multicolumn{1}{ |c| }{coverage} &  0.90 & 0.90 & 0.91 & 0.81 &    \\ 
\wcline{1-6}
\multicolumn{1}{ |c  }{setting 4}  &
\multicolumn{1}{ |c| }{width} & 0.21 & 0.26 & 0.23 & 0.06 &    \\ \cline{2-6}
\multicolumn{1}{ |c  }{$n=100, \ p=95, \ \ \eta=1$}                        &
\multicolumn{1}{ |c| }{coverage} & 0.90 & 0.92 & 0.87 & 0.42 &    \\ \cline{1-6}
\end{tabular}
\label{table:CItable} %
\end{table}
\paragraph{Acknowledgements.} MEL thanks Prof.~Peter J. Bickel for many helpful discussions, and gratefully acknowledges the DOE CSGF under grant  DE-FG02-97ER25308, as well as the NSF-GRFP.

\newpage

  \normalsize

\appendix

\section{Proof of Theorem~\ref{basicprop}}

\proof Due to line~\eqref{decompose} and Lemma 8.8 in B\&F 1981,
\begin{align}\label{biasstep}
d_2^2(\Psi_{\rho}(F_0;c),\Phi_{\rho}(\hat{F};c)) &=  d_2^2\big(\Phi_{\rho}(F_0;c), \Phi_{\rho}(\hat{F};c)\big) + (c\ttop \delta(X))^2.
\end{align}
If $\ve^*\in\R^n$ is a random vector whose entries are drawn i.i.d. from $\hat{F}$, then the definition of $\Phi_{\rho}$ gives the matching relations 
\begin{equation}\label{relation}
\begin{split}
\Phi_{\rho}(F_0;c)&=\mathcal{L}(c\ttop (X\ttop X +\rho)^{-1}X\ttop \ve \big\bracevert X)\\
\Phi_{\rho}(\hat{F};c)&=\mathcal{L}(c\ttop (X\ttop X +\rho)^{-1}X\ttop \ve^*\big\bracevert X).
\end{split}
\end{equation}
To make use of these relations, we apply Lemma 8.9 in B\&F 1981, which implies that if $w\in \R^n$ is a generic deterministic vector, and if $U=(U_1,\dots,U_n)$ and $V=(V_1,\dots,V_n)$ are random vectors with i.i.d. entries, then
$$d_2^2(w\ttop U,w\ttop V) \leq \|w\|_2^2\cdot d_2(U_1,V_1)^2.$$
Therefore, 
\begin{equation}
\begin{split}
d_2^2\big(\Phi_{\rho}(F_0),\Phi_{\rho}(\hat{F})\big) &\leq \|c\ttop (X\ttop X +\rho)^{-1}X\ttop\|_2^2\cdot d_2^2(\ve_1,\ve_1^*)\\[0.2cm]
&=\ts\frac{1}{\sigma^2}\cdot v_{\rho}(X;c)\cdot d_2^2(F_0,\hat{F}).
\end{split}
\end{equation}
Combining this with line~\eqref{biasstep} and dividing through by $v_{\rho}$ proves the claim.\qed

\section{Proof of Theorem~\ref{mspelink}}
\proof
By the triangle inequality, 
\begin{equation}
\begin{split}
d_2^2(\hat{F},F_0)
&\leq 2 \, d_2^2(\hat{F},F_n)+2\,d_2^2(F_n,F_0).
\end{split}
\end{equation}
Let $\tilde{F}_n$ be the (uncentered) empirical distribution of the residuals $\hat{e}$ of $\hat{\beta}$, which places mass $1/n$ at each value $\hat{e}_i$, for $i=1,\dots,n$. The proofs of Lemmas 2.1 and 2.2 in Freedman 1981,  show that
\begin{equation}
\begin{split}
\E \,\big[ d_2(\hat{F},F_n)^2\big \bracevert X \big]&\leq \E \big[\big(\ts\frac{1}{n}\tsum_{i=1}^n \ve_i\big)^2\big] +\E \big[\,d_2(\tilde{F}_n,F_n)^2\big\bracevert X\big]\\[0.2cm]
&\leq \ts\frac{1}{n}\sigma^2+\ts\frac{1}{n}\E\big[\|\hat{e}-\ve\|_2^2\big\bracevert X\big]\\[0.2cm]
&=\ts\frac{1}{n}\sigma^2+\ts\frac{1}{n}\E\big[\|X(\beta-\hatbeta)\|_2^2\big\bracevert X\big],
\end{split}
\end{equation}
where we have used the algebraic identity $\hat{e}-\ve = X(\beta-\hatbeta)$, which holds for any estimator $\hatbeta$. This completes the proof.\qed

\section{Proof of Theorem~\ref{ridgemspe}}\label{proofridgemspe}

\proof

We begin with a simple bias-variance decomposition,
\begin{equation}\label{biasvar}
\begin{split}
\mspe(\hat{\beta}_{\varrho}|X) &= \ts\frac{1}{n}\E\Big[\big\| X\big(\hat{\beta}_{\varrho}-\E\big[\hat{\beta}_{\varrho} | X\big]\big)\big\|_2^2\Big\bracevert X\Big]+\ts\frac{1}{n}\big\|X(\E\big[\hat{\beta}_{\varrho}|X\big]-\beta)\big\|_2^2.
\end{split}
\end{equation}
We will handle the bias and variance terms separately. To consider the bias term, note that
$\E[\hat{\beta}_{\rho}|X]-\beta =B \beta,$
where
$$B=(X\ttop X+\varrho I_{p\times p})^{-1}X\ttop X -I_{p\times p}.$$
Hence,
\begin{equation}
\begin{split}
\ts\frac{1}{n}\|X(\E[\hat{\beta}_{\varrho}|X]-\beta)\|_2^2&= \ts\frac{1}{n}\|XB\beta\|_2^2\\[0.2cm]
&=\beta\ttop B(\ts\frac{1}{n} X\ttop X) B \beta.
\end{split}
\end{equation}
If we let $l_i=\lambda_i(\frac{1}{n}X\ttop X)$, then the eigenvalues of $B(\frac{1}{n}X\ttop X) B$ are of the form $\mu_i:=\frac{l_i(\varrho/n)^2}{(l_i+\varrho/n)^2}$. In particular, it is simple to check\footnote{ Note that if $t\in\R$ and $f(t):=\frac{t(\varrho/n)^2}{(t+\varrho/n)^2}$, then $f$ is maximized at $t=\varrho/n$. Also, if $\theta \leq \nu$, then there at least one $l_i$ that scales at the rate of $\varrho/n$.} that $ \max_i \mu_i \asymp \varrho/n$ whenever $\theta\leq \nu$, and so
\begin{equation}\label{biasrate}
 \beta\ttop B(\ts\frac{1}{n} X\ttop X) B \beta \lesssim \ts\frac{\varrho}{n}\|\beta\|_2^2 = n^{-\theta}\|\beta\|_2^2.
 \end{equation}
Note that this bound is tight, since it is achieved whenever $\beta$ is parallel to the top eigenvector of $B(\frac{1}{n}X\ttop X)B$.

To consider the variance term, note that
$\hat{\beta}_{\varrho}-\E[\hat{\beta}_{\varrho}|X]=(X\ttop X+\varrho I_{p\times p})^{-1}X\ttop \ve,$ and so
\begin{equation}\label{foobar}
\begin{split}
 \ts\frac{1}{n}\E\Big[\big\| X(\hat{\beta}_{\varrho}-\E[\hat{\beta}_{\varrho}|X])\big\|_2^2\Big\bracevert X\Big]&= \ts\frac{1}{n}\tr\Big(\Big( X\ttop X\big(X\ttop X+\varrho I_{p\times p}\big)^{-1}\Big)^2\Big)\\
 &= \ts\frac{1}{n}\sum_{i=1}^{n\wedge p} \big(\ts\frac{l_i}{l_i+\varrho/n}\big)^2.
 \end{split}
 \end{equation}
It is natural to decompose the sum in terms of the index set
\begin{equation}\label{indexset}
\mathcal{I}(n):=\{i\in \{1,\dots,n\wedge p\}: l_i \geq \varrho/n\},
\end{equation}
which satisfies $|\mathcal{I}(n)|\asymp n^{\frac{\theta}{\nu}}$. We will bound the variance term in two complementary cases; either $\nu>1/2$ or $\nu\leq 1/2$. First assume $\nu>1/2$. Then,

\begin{align}\label{sum}
\ts\frac{1}{n}\displaystyle \sum_{i=1}^p \ts\big(\frac{ l_i}{ l_i+\varrho/n}\big)^2&=\ts\frac{1}{n}\displaystyle \sum_{i\in \mathcal{I}(n)} \ts\big(\frac{ l_i}{ l_i+\varrho/n}\big)^2
+\ts\frac{1}{n}\displaystyle \sum_{i\not\in \mathcal{I}(n)} \ts\big(\frac{ l_i}{ l_i+\varrho/n}\big)^2\\[0.2cm]
&\lesssim \ts\frac{1}{n}\,|\mathcal{I}(n)|+\ts\frac{1}{n}\displaystyle \int_{|\calI(n)|}^{n\wedge p}\ts\frac{x^{-2\nu}}{(\varrho/n)^2}dx\\[0.2cm]
&\lesssim n^{-1}\Big( n^{\frac{\theta}{\nu}}+n^{2\theta}\cdot  (|\calI(n)|)^{1-2\nu}\Big) \ \ \ \text{ using } \nu>\ts\frac{1}{2}\label{alphahalf}\\[0.2cm]
&\asymp n^{-1}\Big( n^{\frac{\theta}{\nu}}+n^{2\theta}\cdot (n^{\frac{\theta}{\nu}})^{(1-2\nu)}\Big)\\
&= 2n^{\frac{\theta-\nu}{\nu}}\label{last}.
\end{align}
To see that this upper bound is tight, note that in line~\eqref{sum}, we can use the term-wise lower bounds 
\begin{equation}
\big(\ts\frac{l_i}{l_i+\varrho/n}\big)^2\geq
\begin{cases}
&\ts\frac{1}{4} \, \ \ \ \ \ \ \ \ \ \, \,  \text{ if } i\in \calI(n)\\
& \ts\frac{1}{4}\frac{l_i^2}{(\varrho/n)^2} \, \ \ \text{ if } i\not\in \calI(n),
\end{cases}
\end{equation}
and then apply an integral approximation from below (which leads to the same rate).
Combining the bias and variance pieces, we have shown that 
$$\ts\frac{1}{n}\E\|X(\hat{\beta}_{\varrho}-\beta)\|_2^2 \lesssim n^{\frac{\theta-\nu}{\nu}}+n^{-\theta} \text{ \ \ \ \ \ if  \ \ \ } \nu>\frac{1}{2}.$$
 The bound is optimized when the two terms on the right side have the same rate, which leads to the choice $\theta=\frac{\nu}{\nu+1}$. 

In the case where $\nu\in(0,\ts\frac{1}{2})$, the calculation proceeds in the same way up to line~\eqref{alphahalf}, where we obtain the bound
\begin{align}\label{sumagain}
\ts\frac{1}{n}\displaystyle \sum_{i=1}^{n\wedge p} \ts\big(\frac{ l_i}{ l_i+\varrho/n}\big)^2 
&\lesssim n^{-1}\Big( n^{\frac{\theta}{\nu}}+n^{2\theta}\cdot  n^{1-2\nu}\Big)\\
&= n^{\frac{\theta-\nu}{\nu}}+n^{2(\theta-\nu)}\label{alphabigger}.
\end{align}
This bound is also tight due to the same reasoning as above. Note that in order for the bound to tend to 0 as $n\to\infty$, we must choose $\theta<\nu$. Furthermore, since we are working under the assumption $\nu\in (0,\ts\frac{1}{2})$, it follows that the right side of line~\eqref{alphabigger} has rate equal to $n^{2(\theta-\nu)}$. Combining the rates for the bias and variance shows that 
$$\ts\frac{1}{n}\E\|X(\hat{\beta}_{\varrho}-\beta)\|_2^2 \lesssim n^{2(\theta-\nu)}+n^{-\theta}
\text{ \ \ \ \ \ if  \ \ \ } \nu\in(0,\frac{1}{2}).
$$
 The bound is optimized when the two terms on the right side have the same rate, which leads to the choice $\theta=\frac{2\nu}{3}$. 
\qed

\section{Proof of Lemma~\ref{propGaussian}}\label{appPropGaussian}

The proof is split up into three pieces, corresponding to parts (i), (ii), and (iii) in the statement of the result.

\subsection{The bias inequality~\eqref{biasineq} }
We prove inequality~\eqref{biasineq} by combining Lemmas~\ref{biasrep} and~\ref{expectedbias} below.

\begin{lemma}\label{biasrep}
Assume the conditions of Lemma~\ref{propGaussian}. For each $i\in \{1,\dots, n\}$, there are independent random vectors $u_i(X), w(X)\in \R^p$ such that the random variable $X_i\ttop\delta(X)$ can be represented algebraically as
$$b_{\rho}(X;X_i)=X_i\ttop\delta(X) = u_i(X)\ttop w(X).$$
Here, the vectors $u_i(X)$ can be represented in law as
\begin{equation}
u_i(X)  \overset{\mathcal{L}}{=} \ts\frac{1}{\|z\|_2} \Pi_p(z),
\end{equation}
where $z\in\R^n$ is a standard Gaussian vector, and $\Pi_p(z):=(z_1,\dots,z_p)$. Also, the vector $w(X)$ satisfies the bound $\|w(X)\|_2^2 \leq \ts\frac{\rho}{4}\|\beta\|_2^2$ almost surely.
\end{lemma}
\proof To fix notation, we write $X\ttop = \Sigma^{1/2} Z\ttop$ where $Z\ttop\in \R^{p\times n}$ is a standard Gaussian matrix. Recall that $\delta(X) = B\beta,$
where 
$$B= I_{p\times p} - (X\ttop X+\rho I_{p\times p})^{-1}X\ttop X.$$
Let $Z= HLG\ttop$ be a ``signed s.v.d.''~for $Z$, as defined in Appendix~\ref{sec:gaussian}, where $H\in\R^{n\times p}$, $L\in\R^{p\times p}$, and $G\in\R^{p\times p}$.
Now define $u_i(X)$ and $w(X)$ according to
\begin{align}
X_i\ttop \delta(X) = e_i\ttop X B\beta=e_i\ttop Z \Sigma^{1/2}B\beta=\underbrace{e_i\ttop H}_{=:\,u_i(X)\ttop}\!  \underbrace{LG\ttop \Sigma^{1/2} B\beta}_{=:\ts w(X)}.
\end{align}
 From Lemma~\ref{svdlemma} in Appendix~\ref{sec:gaussian}, the rows $e_i\ttop H$ can be represented in distribution as $\frac{1}{\|z\|_2}\Pi_p(z)$. The same lemma also shows that the three matrices $H$, $L$, and $G$ are independent.

Hence, to show that $u_i(X)$ and $w(X)$ are independent, it suffices to show that $w(X)$ is a function only of $G$ and $L$. In turn, it is enough to show that $B$ is a function only of $G$ and $L$. But this is simple, because $B$ is a function only of the matrix $X\ttop X$, which may be written as
\begin{equation}\label{factor}
X\ttop X =\Sigma^{1/2}Z\ttop Z \Sigma^{1/2} = \Sigma^{1/2} GL^2G\ttop \Sigma^{1/2}.
\end{equation}
It remains to show that $\|w(X)\|_2^2 \leq \ts\frac{\rho}{4}\|\beta\|_2^2$ almost surely.  Combining the definition of $w(X)$ with line~\eqref{factor} gives
\begin{align}
\|w(X)\|_2^2 =  \beta \ttop \big(B  X\ttop X B \big)\beta.
\end{align}
The eigenvalues of $ BX\ttop X B$ are of the form $\mu_i:= n\frac{(\rho/n)^2 l_i}{(l_i+\rho/n)^2}$ where $l_i=\lambda_i(\frac{1}{n}X\ttop X)$, and it is simple to check that the inequality $\max_i \mu_i \leq \ts\frac{\rho}{4}$ holds for every realization of $X$. 
\qed\\

Before proceeding to the second portion of the proof of inequality~\eqref{biasineq}, we record some well-known tail bounds for Gaussian quadratic forms due to Laurent and Massart~\cite{LaurentMassart}, which will be useful at various points later on.

\begin{lemma}[Laurent \& Massart, 2001]\label{laurentmassart}
Let $A\in \R^{n \times n}$ be a fixed symmetric matrix, and let $z\in \R^n$ be a standard Gaussian vector. Then, for every $t>0$,
\begin{equation}\label{lm1}
\P\Big[z\ttop A z \geq \tr(A)+2\mnorm{A}_{F}\sqrt{t}+2\mnorm{A}_{\textup{op}} t \Big]\leq \exp(-t)
\end{equation}
and
\begin{equation}\label{lm2}
\P\Big[ z\ttop A z \leq \tr(A)-2\mnorm{A}_F\sqrt{t}\Big]\leq \exp(-t).
\end{equation}
\end{lemma}
~\\
The next lemma completes the proof of inequality ~\eqref{biasineq}.
\begin{lemma}\label{expectedbias}
Assume the conditions of Lemma~\ref{propGaussian}, and let $\tau> 0$ be a constant. Then for every $n\geq 1$, the following event holds with probability at least $1-n^{-\tau}-ne^{-n/16}$, 
\begin{equation}
\max_{1\leq i\leq n} b_{\rho}^2(X;X_i)\ \leq 5\|\beta\|_2^2\cdot n^{-\gamma} \cdot  (\tau+1)\log(n+2).
\end{equation}
\end{lemma}

\proof Applying the representation for $b_{\rho}(X;X_i)$ given in Lemma~\ref{biasrep}, there is a standard Gaussian vector $z\in \R^n$, such that $u_i(X)\de \Pi_p(z)/\|z\|_2$. Consequently,
\begin{align}
b_{\rho}^2(X;X_i) 
&\de\ts\frac{1}{\|z\|_2^2} \cdot \Pi_p(z) \Big(w(X)w(X)\ttop \Big) \Pi_p(z),
\end{align}
where we may take $z$ and $w(X)$ to be independent by the same lemma. Using Lemma~\ref{laurentmassart} on Gaussian quadratic forms, as well as the fact that $\|w(X)\|_2^2 \leq \frac{\rho}{4} \|\beta\|_2^2$ almost surely,  we have for all $t>0$,
\begin{equation}\label{numbound}
\P\Bigg[\Pi_p(z) \ttop \Big(w(X)w(X)\ttop \Big) \Pi_p(z) \geq \ts\frac{\rho}{4} \|\beta\|_2^2 \big(1+2\sqrt{t}+2t)\Bigg\bracevert w(X)\Bigg]\leq \exp(-t).
\end{equation}
The same lemma also implies that for all $t'\in (0,\frac{1}{4})$,
\begin{equation}\label{denombound}
\P\Big[ \ts\frac{1}{\|z\|_2^2} \geq \ts\frac{1}{(1-2\sqrt{t'})n}\Big] \leq \exp(-nt').
\end{equation}

Now, we combine the bounds by integrating out $w(X)$ in line~\eqref{numbound} and choosing $t'=1/16$ in line~\eqref{denombound}. Taking a union bound, we conclude that for any $t>0$, and any fixed $i=1,\dots,n$,
\begin{equation}
\P\Bigg[ b_{\rho}^2(X;X_i)\ \leq \ \ts\frac{\rho}{n} \cdot \|\beta\|_2^2\cdot \ts \frac{1}{2} \big(1+2\sqrt{t}+2t\big)\Bigg]\geq 1-e^{-t}-e^{-n/16}.
\end{equation}
Finally, another union bound shows that the maximum of the $b_{\rho}(X,X_i\ttop)$ satisfies

\begin{equation}
\P\Bigg[ \max_{1\leq i\leq n} b_{\rho}^2(X;X_i)\ \leq \ \ts\frac{\rho}{n} \cdot \|\beta\|_2^2\cdot \ts\frac{1}{2}(1+2\sqrt{t}+2t)\Bigg]\geq 1-e^{-t+\log(n)}-ne^{-n/16},
\end{equation}
which implies the stated result after choosing $t=(\tau+1) \log(n+2)$, and noting that since $t\geq 1$, we have $\ts\frac{1}{2}(1+2\sqrt{t}+2t)\leq 5t  = 5 (\tau+1) \log(n+2)$, as well as $e^{-t+\log(n)}\leq e^{-\tau \log(n+2)}\leq n^{-\tau}$ for every $n\geq1$.
\qed

\subsection{The variance inequality~\eqref{varianceineq} }\label{variancesubsec}
The following ``representation lemma'' will serve as the basis for controlling the variance $v_{\rho}(X;X_i) = \sigma^2 \|X_i\ttop (X\ttop X+\rho I_{p\times p})^{-1}X\ttop\|_2^2$.

\begin{lemma}\label{variancerep}
Assume the conditions of Lemma~\ref{propGaussian}. For each $i\in \{1,\dots,n\}$, there is a random vector $v_i(X)\in \R^p$ and a random matrix $M(X)\in \R^{p\times p}$ that are independent and satisfy the algebraic relation
$$ \|X_i\ttop (X\ttop X+\rho I_{p\times p})^{-1}X\ttop\|_2^2 = v_i(X)\ttop M(X) v_i(X).$$
Here, the vector $v_i(X)$ can be represented in law as
\begin{equation}
v_i(X)\overset{\mathcal{L}}{=} \ts\frac{1}{\|z\|_2} \Pi_p(z),
\end{equation}
where $z\in\R^n$ is a standard Gaussian vector and $\Pi_p(z)=(z_1,\dots,z_p)$. Also,
the matrix $M(X)$ satisfies the algebraic relation
\begin{equation}\label{tracerelation}
\tr(M(X)) = \|X (X\ttop X+\rho I_{p\times p})^{-1}X\ttop\|_F^2.
\end{equation}
\end{lemma}
An explicit formula for $M(X)$ is given below.
\proof Define the matrix $A:=(X\ttop X+\rho I_{p\times p})^{-1}X\ttop X(X\ttop X+\rho I_{p\times p})^{-1}$. Then,
\begin{align}
\|X_i\ttop (X\ttop X+\rho I_{p\times p})^{-1}X\ttop\|_2^2 
&=e_i\ttop X AX\ttop e_i.
\end{align}
Using the notation in the proof of the previous lemma, let $X = Z\Sigma^{1/2}$ where $Z\in \R^{n\times p}$ is a standard Gaussian random matrix. Furthermore, let $Z = H L G\ttop$ be a signed s.v.d. for $Z$, as defined in Appendix~\ref{sec:gaussian}.
Then, we define $v_i(X)$ and $M(X)$ according to
\begin{align}
e_i\ttop XAX\ttop e_i
&=\underbrace{e_i\ttop H}_{=:v_i(X)\ttop} \underbrace{LG\ttop \Sigma^{1/2} A \Sigma^{1/2}GL\ttop}_{=:M(X)}H\ttop e_i.
\end{align}
Some algebra shows that $M(X)$ satisfies the relation~\eqref{tracerelation}. As in the proof of Lemma~\ref{biasrep}, the argument is completed using two properties of the signed s.v.d.~of a standard Gaussian matrix: The rows of $H$ can be represented as $\Pi_p(z)/\|z\|_2$ where $z\in\R^n$ is a standard Gaussian vector, and the matrices $H$, $L$, and $G\ttop$ are independent. (See Lemma~\ref{svdlemma} in Appendix~\ref{sec:gaussian}.) To show that $v_i(X)$ and $M(X)$ are independent, first note that $v_i(X)$ only depends on $H$. Also, it is simple to check that $M(X)$ only depends on $G$ and $L$, because $A$ is a function only of $X\ttop X= \Sigma^{1/2}GL^2G\ttop \Sigma^{1/2}$.
\qed

\subsubsection{Concentration of the variance and bounds on its expected value}
Due to Lemma~\ref{variancerep}, for each $i=1,\dots,n$, we have the representation
\begin{equation}
v_{\rho}(X;X_i)
\de \ts\frac{1}{\|z\|_2^2} \Pi_p(z)\ttop M(X) \Pi_p(z),
\end{equation}
where $z\in \R^n$ is a standard Gaussian vector, independent of $M(X)$. Conditionally on $M(X)$, the quadratic form $\Pi_p(z)\ttop M(X) \Pi_p(z)$ concentrates around $\tr(M(X))$ by Lemma~\ref{laurentmassart}. The same lemma also implies that $\|z\|_2^2$ concentrates around $n$. In the next three subsections, we will show that $\sqrt{\tr(M(X))}$ concentrates around its expected value, and obtain upper and lower bounds on the expected value. We will need two-sided bounds in preparation for Theorem~\ref{lasthm}.

\subsubsection{Concentration of $\sqrt{\tr(M(X))}$}

\begin{lemma}\label{sqconc}
Assume the conditions of Lemma~\ref{propGaussian}. Then for every $t>0$, and every $n\geq 1$, 
\begin{equation}
\P\Bigg[ \Big|\sqrt{\tr(M(X))}- \E\sqrt{\tr(M(X))}\Big| \geq t\Bigg] \leq 2\exp(-\ts\frac{64}{54}\ts \frac{n^{1-\gamma}t^2}{ \|\Sigma\|_{\textup{op}}}).
\end{equation}
\end{lemma}
\proof We will show that $\sqrt{\tr(M(X))}$ is a Lipschitz function of a standard Gaussian matrix.
Define the function $g_{\rho}:\R_+\to [0,1]$ by $g_{\rho}(s)=\frac{s^2}{s^2+(\rho/n)}$, which satisfies the Lipschitz condition
$$|g_{\rho}(s)-g_{\rho}(s')|\leq\mathfrak{L}_n|s-s'|,$$
for all $s,s'\geq 0$, where $\mathfrak{L}_n:=\ts\frac{3\sqrt{3}}{8} \ts\frac{1}{\sqrt{\rho/n}}$.

If $\sigma(A)=(\sigma_1(A),\dots,\sigma_k(A))$ denotes the vector of singular values of a rank $k$ matrix $A$, then we define $g_{\rho}$ to act on $\sigma(A)$ component-wise, i.e. $g_{\rho}(\sigma(A)) = (g_{\rho}(\sigma_1(A)),\dots,g_{\rho}(\sigma_k(A)))$. Recall from Lemma~\ref{variancerep} that
\begin{equation}
\begin{split}
\sqrt{\tr(M(X))}&= \|X(X\ttop X+\rho I_{p\times p})^{-1}X\ttop \|_F\\
\end{split}
\end{equation}
and note that the  $i$th singular value of the matrix $X(X\ttop X+\rho I_{p\times p})^{-1}X\ttop$ is given by $g_{\rho}(\sigma_i(\ts\frac{1}{\sqrt{n}}X))$. Viewing the Frobenius norm of a matrix as the $\ell_2$ norm of its singular values, we have
\begin{equation}
\begin{split}
\sqrt{\tr(M(X))} &= \ts\|g_{\rho}(\sigma(\ts\frac{1}{\sqrt{n}}X))\|_2.
\end{split}
\end{equation}
Write $X\ttop =\Sigma^{1/2}Z\ttop $ for a standard Gaussian matrix $Z\in \R^{n\times p}$,
 and let $f:\R^{n\times p}\to\R$ be defined according to
$$f(Z) := \sqrt{\tr(M(X))}.$$
We claim that $f$ is Lipschitz with respect to the Frobenius norm.
 Let $W\ttop \in \R^{p\times n}$ be a generic matrix, and put $A=\ts\frac{1}{\sqrt{n}}\Sigma^{1/2}Z\ttop$ and $B=\ts\frac{1}{\sqrt{n}}\Sigma^{1/2}W\ttop$.
 Then,
\begin{align}
|f(Z)-f(W)| &= \Big| \|g_{\rho}(\sigma(A))\|_2 - \|g_{\rho}(\sigma(B))\|_2\Big|\\[0.2cm]
&\leq \big\|g_{\rho}(\sigma(A)) -g_{\rho}(\sigma(B))\big\|_2\\[0.2cm]
&\leq \mathfrak{L}_n \|\sigma(A)-\sigma(B)\|_2\\[0.2cm]
&\leq  \mathfrak{L}_n\|A-B\|_F   \ \ \  \ \ \ \ \text{(Weilandt-Hoffman)}\\[0.2cm]
&= \ts  \mathfrak{L}_n\big\|\ts\frac{1}{\sqrt{n}}\Sigma^{1/2}\big(Z\ttop-W\ttop \big)\big\|_F  \\[0.2cm]
&\leq \ts\frac{\mathfrak{L}_n}{\sqrt{n}}\sqrt{\|\Sigma\|\op}\cdot \big\|Z\ttop-W\ttop \big\|_F, 
\end{align}
where we have used a version of the Weilandt-Hoffman inequality for singular values~\cite[p.186]{horntopics}, as well as the inequality $\|M_1M_2\|_F\leq \|M_1\|\op\|M_2\|_F$, which holds for any square matrix $M_1$ that is compatible with $M_2$. (See Lemma~\ref{frobineq} in Appendix~\ref{appmisc}.)
The statement of the lemma now follows from the Gaussian concentration inequality. (See Lemma~\ref{GaussianConc} in Appendix~\ref{appmisc}).\qed

\subsubsection{Upper bound on $\E\sqrt{\tr(M(X))}$ }

\begin{lemma}\label{traceupper} Assume the conditions of Lemma~\ref{propGaussian}. Then, the matrix $M(X)$ satisfies
\begin{equation}
\E\sqrt{\tr(M(X))}\lesssim
\begin{cases}
   & n^{(\gamma-\eta)+\frac{1}{2}} \text{\ \ \ \  \ if } \ \eta \in (0, \ts\frac{1}{2})\\
   & n^{\frac{\gamma}{2\eta}} \text{\   \ \ \  \  \ \ \  \ \ \,  if } \ \eta>\ts\frac{1}{2}.
   \end{cases}
\end{equation}
\end{lemma}

\proof By Jensen's inequality, it is enough to bound $\sqrt{\E[\tr(M(X))]}$ from above. Define the univariate function $\psi: \R_+\to \R_+$ by $\psi(s):=\frac{s}{(\sqrt{s}+\rho/n)^2}$, and observe that
\begin{equation}\label{tracecalc}
\begin{split}
\tr(M(X)) &= \tr\Big(\Big((X\ttop X+\rho I_{p\times p})^{-1}X\ttop X\Big)^2\Big)\\
&=\sum_{i=1}^p \ts\frac{\lambda_i^2(\hat{\Sigma})}{(\lambda_i(\hat{\Sigma})+\rho/n)^2}\\
&=\sum_{i=1}^p \psi(\lambda_i(\hat{\Sigma}^2))\\
&=\tr\big(\psi(\hat{\Sigma}^2)\big).
\end{split}
\end{equation}
Here where we use the ``operator calculus'' notation $\psi(A)=U\psi(D)U\ttop$ where $A$ is a symmetric matrix with spectral decomposition $A=UDU\ttop$, and $\psi(D)$ is the diagonal matrix whose $i$th diagonal entry is $\psi(D_{ii})$.
 It is simple to check that $\psi$ is a concave, and so $\tr(\psi(\hat{\Sigma}^2))$ is a concave matrix functional of $\hat{\Sigma}^2$ by Lemma~\ref{convexfunctional} in Section~\ref{convex} of Appendix A. Therefore, Jensen's inequality implies
\begin{equation}\label{jensensq}
\begin{split}
\E[\tr(M(X))] &\leq \tr(\psi(\E[\hat{\Sigma}^2]))\\[0.2cm]
&=\ts\sum_{i=1}^p \psi(\lambda_i(\mathfrak{S})),
\end{split}
\end{equation}
where we define the matrix $\mathfrak{S}:=\E[\hat{\Sigma}^2]$. Since $X$ is Gaussian, $\hat{\Sigma}$ is a Wishart matrix up to scaling, and so Lemma~\ref{wishartsquared} in Appendix~\ref{appmisc} shows that this expectation may be evaluated exactly as
\begin{equation}
\mathfrak{S}=(1+\ts\frac{1}{n})\Sigma^2+\ts\frac{\tr(\Sigma)}{n}\Sigma.
\end{equation}
We will now use this relation to apply an integral approximation to the right side of line~\eqref{jensensq}.
Clearly, the eigenvalues of $\mathfrak{S}$ are given by
\begin{equation}
\begin{split}
\lambda_i(\mathfrak{S})&=(1+\ts\frac{1}{n})\lambda_i^2(\Sigma)+\ts\frac{\tr(\Sigma)}{n}\lambda_i(\Sigma)\\[0.2cm]
&\asymp i^{-2\eta}+\ts\frac{\tr(\Sigma)}{n}i^{-\eta}.
\end{split}
\end{equation}
Let $r\in (0,1)$ be a constant to be specified later. On the set of indices $1\leq i\leq \lceil n^r\rceil$ we  use the bound $\psi(\lambda_i(\mathfrak{S}))\leq 1$, and on the set of indices $i >\lceil n^r\rceil$ we use the bound $\psi(\lambda_i(\mathfrak{S}))\leq \frac{1}{(\rho/n)^2} \lambda_i(\mathfrak{S})$. Recalling the assumption $\rho/n=n^{-\gamma}$, we may decompose the inequality~\eqref{jensensq} as\footnote{Note that if $p/n\asymp 1$, it is possible that $n^r>p$ for small values of $n$. Since we want $n^r\leq p$ for the integral in line~\eqref{makesense}, Lemma~\ref{propGaussian} is stated for ``all large $n$''.}
\begin{align}
\E[\tr(M(X))] &\leq \ \sum_{i=1}^{\lceil n^r\rceil} \psi(\lambda_i(\mathfrak{S}))+ \sum_{i=\lceil n^r\rceil+1}^{p} \psi(\lambda_i(\mathfrak{S}))\\[0.2cm]\label{makesense}
&\lesssim n^r+ n^{2\gamma}\int_{n^r}^p \big( x^{-2\eta}+\ts\frac{\tr(\Sigma)}{n} x^{-\eta}\big)dx\\[0.3cm]
&=: n^r+n^{2\gamma} h_n(\eta,r).\label{splitup}
\end{align}
where the function $h_n$ is defined in the last line.
The bound is optimized when the two terms on the right are of the same order; i.e.~when $r$ solves the rate equation
\begin{equation}\label{rateeqn}
n^r\asymp n^{2\gamma} h_n(\eta,r).
\end{equation}
Noting that
\begin{equation}
\ts\frac{\tr(\Sigma)}{n}\asymp
\begin{cases}
n^{-\eta} & \text{ if } \eta\in (0,1)\\
n^{-1} & \text{ if  } \eta>1,
\end{cases}
\end{equation}
the quantity $h_n(\eta,r)$ may be computed directly as
\begin{equation}
h_n(\eta,r)\asymp
\begin{cases}
   n^{1-2\eta}  \text{\ \ \ \ \ \  \ \ \ \ \ \ \ \  \ \ \ \ \ \ \ \  \ \ \ \ \ \ \, if } \ \eta \in (0, \ts\frac{1}{2}),\\
    n^{r(1-2\eta)} \text{\ \ \ \  \ \ \ \ \ \ \ \ \ \ \ \ \ \ \ \ \ \ \  \ \,  if } \ \eta \in (\ts\frac{1}{2}, 1),\\
    n^{r(1-2\eta)}+n^{r(1-\eta)-1} \text{\ \ \ \ if } \ \eta>1.
   \end{cases}
\end{equation}
If we let $r=r_*(\eta,\gamma)$ denote the solution of the rate equation~\eqref{rateeqn}, then some calculation shows that under the assumption $\gamma\in (0,1)$,
\begin{equation}
r_*(\eta,\gamma)=
\begin{cases}
   2(\gamma-\eta)+1 \text{\ \ \ \ \  if } \ \eta \in (0, \ts\frac{1}{2}),\\
    \ts\frac{\gamma}{\eta} \text{\   \ \ \  \ \ \  \ \ \ \ \ \ \ \ \ \ \ \ \  \ \ \   if } \ \eta>\ts\frac{1}{2}.
   \end{cases}
\end{equation}
When $\eta\in(0,1)$ this is straightforward. To show the details for $\eta>1$, note that the rate equation~\eqref{rateeqn} may be written as
\begin{equation}\label{firstform}
n^r \asymp n^{2\gamma +r(1-2\eta)}+n^{2\gamma+r(1-\eta)-1},
\end{equation}
which is the same as
\begin{equation}\label{secondform}
1 \asymp n^{2(\gamma-\eta r)}+n^{2\gamma-\eta r-1}.
\end{equation}
In order for both terms on the right to be $\mathcal{O}(1)$, the number $r$ must satisfy the constraints
\begin{align}
r &\geq \ts\frac{\gamma}{\eta}, \\
r &\geq \ts\frac{\gamma}{\eta}+\ts\frac{\gamma-1}{\eta}.
\end{align}
Since Lemma~\ref{propGaussian} assumes $\gamma\in (0,1)$, only the first constraint matters. Furthermore, when $r\geq \frac{\gamma}{\eta}$, the second term in line~\eqref{secondform} is $o(1)$, and we are reduced to choosing $r$ so that $1\asymp n^{2(\gamma-\eta r)}$, which gives $r=r_*(\eta,\gamma)=\frac{\gamma}{\eta}$. Substituting this value into line~\eqref{splitup} completes the proof. (Note from the discussion preceding line~\eqref{splitup} that $r$ must lie in the interval $(0,1)$, and this requires $\gamma/\eta<1$, which explains the assumption $\gamma<\min\{\eta,1\}$ in Lemma~\ref{propGaussian}.)
\subsubsection{Lower bound on $\E\sqrt{\tr(M(X))}$ }

\begin{lemma}\label{tracelower}
Assume the conditions of Lemma~\ref{propGaussian}. Then, the matrix $M(X)$ satisfies
\begin{equation}
\E\sqrt{\tr(M(X))} \, \gtrsim \, n^{\frac{\gamma}{2\eta}}.
\end{equation}¥
\end{lemma}

\proof 

The variable $\sqrt{\tr(M(X))}$ may be written as $\|X\ttop X(X\ttop X+\rho I_{p\times p})^{-1}\|_F.$
Since the Frobenius norm is a convex matrix functional, Jensen's inequality implies
\begin{equation}
\begin{split}
\E\sqrt{\tr(M(X))}&\geq \Big\|\E\Big[X\ttop X(X\ttop X+\rho I_{p\times p})^{-1}\Big]\Big\|_F\\[0.2cm]
&=\Big\|\E\Big[\big(I_{p\times p}+\ts\frac{\rho}{n}\hat{\Sigma}^{-1}\big)^{-1}\Big]\Big\|_F,
\end{split}
\end{equation}
where the last step follows algebraically with $\hat{\Sigma}:=\ts\frac{1}{n}X\ttop X$. If we define the univariate function $f:\R_+\to\R_+$ by
$f(s)=(1+\frac{\rho}{n}s)^{-1}$, then last inequality is the same as
\begin{equation}
\E\sqrt{\tr(M(X))}\geq\big\| \E\big[f\big(\hat{\Sigma}^{-1}\big)\big] \big\|_F.
\end{equation}
It is a basic fact that $f$ is operator convex on the domain of positive semidefinite matrices~\cite[p.117]{bhatia}. This yields an operator version of Jensen's inequality with respect to the Loewner ordering (Lemma~\ref{opjensen} in Appendix~\ref{appmisc}):
\begin{equation}
\E\big[f\big(\hat{\Sigma}^{-1}\big)\big] \succeq f\big(\E\big[\hat{\Sigma}^{-1}\big]\big).
\end{equation}
Furthermore, if two matrices satisfy $A\succeq B\succeq 0$, then $\|A\|_F \geq \|B\|_F$~\cite[Corollary 7.7.4]{hornjohnson}. Using this fact, as well as the formula for the expected inverse of a Wishart matrix~\cite[p. 97]{muirhead}, we obtain
\begin{equation}
\begin{split}
\E\sqrt{\tr(M(X))}&\geq \big\| f\big(\E\big[\hat{\Sigma}^{-1}\big]\big)\big\|_F\\[0.2cm]
&=\big\|f\big(\ts\frac{n}{n-p-1}\Sigma^{-1}\big)\big\|_F\\[0.2cm]
&=\Bigg(\sum_{i=1}^p \frac{1}{\big(1+\frac{\rho}{n}\cdot \frac{n}{n-p-1} \lambda_i(\Sigma^{-1})\big)^2}\Bigg)^{1/2}\\[0.2cm]
&=\Bigg(\sum_{i=1}^p \frac{\lambda_i^2(\Sigma)}{\big(\lambda_i(\Sigma)+\frac{\rho}{n}\cdot \frac{n}{n-p-1} \big)^2}\Bigg)^{1/2}. 
\end{split}
\end{equation}
Define the index set $J=\big\{i\in \{1,\dots,p\}: \lambda_i(\Sigma) \geq \frac{\rho}{n}\frac{n}{n-p-1}\big\}$. For any $i\in J$, the $i$th summand in the previous line is at least $1/4$. Also,  assumption~{\textbf{A\ref{Adim}} that $p/n$ is bounded strictly between 0 and 1, as well as the decay condition  on the $\lambda_i(\Sigma)$, imply that  $|J|\asymp n^{\gamma/\eta}$, which completes the proof.
\qed
\subsubsection{Putting the variance pieces together}
Combining Lemmas~\ref{sqconc},~\ref{traceupper}, and~\ref{tracelower} with the Gaussian concentration inequality  (Lemma~\ref{GaussianConc} in Section~\ref{sec:gaussian} of Appendix~A) immediately gives the following result. (We choose $t$ to be proportional to the relevant bound on $\E[\sqrt{\tr(M(X))}]$ in the Gaussian concentration inequality.)

\begin{lemma}\label{traceconc}
Assume the conditions of Lemma~\ref{propGaussian} and let $\tr(M(X))$ be as in line~\eqref{tracerelation}. Then, there are absolute constants $\kappa_1,\kappa_2,\dots,\kappa_6>0$ such that the following upper-tail bounds hold for all large $n$,
\begin{equation}
\P\Big[ \tr(M(X)) \geq \kappa_1 n^{2(\gamma-\eta)+1} \Big] \leq \exp(-\kappa_2 n^{2(1-\eta)+\gamma}), \text{ \ \ \  if \ } \eta \in (0,\ts\frac{1}{2}),
\end{equation}
and
\begin{equation}
\P\Big[ \tr(M(X)) \geq \kappa_3 n^{\gamma/\eta} \Big] \leq \exp(-\kappa_4 n^{1+\frac{\gamma(1-\eta)}{\eta}}) , \text{ \ \ \ \ if \ } \eta >\ts\frac{1}{2},
\end{equation}
and the following lower-tail bound holds for all large $n$,
\begin{equation}
\P\Big[ \tr(M(X)) \leq \kappa_5 n^{\gamma/\eta} \Big] \leq \exp(-\kappa_6 n^{1+\frac{\gamma(1-\eta)}{\eta}}) , \text{ \ \ \ \ if \ } \eta >0.
\end{equation}
\end{lemma}
\paragraph{Remarks.}
Note that in order for the last two probabilities to be small for large values of $\eta>0$, it is necessary that $\gamma<1$, as assumed in Lemma~\ref{propGaussian}. 
The next result completes the assembly of the results in this Subsection~\ref{variancesubsec}. Although the first two bounds in Lemma~\ref{varconc} are not necessary for the statement of Theorem~\ref{lasthm}, they show that the variance $v_{\rho}(X;X_i)$ tends 0 as $n\to\infty$ when $\gamma<\eta$, as assumed in Theorem~\ref{lasthm}. In other words, we imposed the assumption $\gamma<\eta$ so that confidence intervals based on $\Phi_{\rho}(\hat{F}_{\varrho};X_i)$ have width that tends to 0 asymptotically.
\begin{lemma}\label{varconc}
Assume the conditions of Theorem~\ref{lasthm} and let $\tr(M(X))$ be as in line~\eqref{tracerelation}. Assume $\gamma< \min\{\eta,1\}$. Then, there are absolute constants $k_1,k_2,\dots,k_6>0$ such that the following upper-tail bounds hold for all large $n$,
\begin{equation}
\P\Big[ \max_{1\leq i\leq n} v_{\rho}(X;X_i)  \leq  k_1 n^{2(\gamma-\eta)} \Big] \geq 1-4n\exp(-k_2 n^{\frac{\gamma}{\eta}} ) , \text{ \ \ \ \ if \ } \eta\in(0,\ts\frac{1}{2})
\end{equation}
and
\begin{equation}
\P\Big[ \max_{1\leq i\leq n} v_{\rho}(X;X_i)  \leq  k_3 n^{\frac{\gamma}{\eta}-1}  \Big] \geq 1- 4n\exp(-k_4 n^{\frac{\gamma}{\eta}}) , \text{ \ \ \ \ if \ } \eta >\ts\frac{1}{2},
\end{equation}
and
\begin{equation}\label{recip}
\P\Big[ \max_{1\leq i\leq n}\ts \frac{1}{v_{\rho}(X;X_i)}  \leq  k_5  n^{1-\frac{\gamma}{\eta}}\Big] \geq 1- 4n\exp(-k_6 n^{\frac{\gamma}{\eta}}) , \text{ \ \ \ \ if \ } \eta > 0.%
\end{equation}

\end{lemma}

\proof We only prove the last inequality~\eqref{recip}, since the other two inequalities are proven in a similar way.
By Lemma~\ref{variancerep}, we have 
\begin{equation}\label{quadform}
v_{\rho}(X;X_i) \de \ts\frac{1}{\|z\|_2^2} \Pi_p(z)\ttop M(X) \Pi_p(z)
\end{equation}
 where $z\sim N(0,I_{p\times p})$ and $z \independent M(X)$. To apply the lower-tail bound for Gaussian quadratic forms, note that H\"older's inequality implies $\|M(X)\|_F\leq \sqrt{\tr(M(X))}$ since $\|M(X)\|\op\leq 1$ almost surely. Therefore, letting $t=t'\tr(M(X))$ with $t'\in (0,1)$ in inequality~\eqref{lm2} gives
 \begin{equation}
\P\Big[ \Pi_p(z)\ttop M(X) \Pi_p(z) \geq (1-2\sqrt{t'})\tr(M(X))\Big\bracevert M(X) \Big] \geq 1-\exp\big(-t' \tr(M(X))\big)
\end{equation}
Next, observe that inequality~\eqref{lm2} with $t=t'\cdot n$ for $t'\in (0,1)$ gives,
\begin{equation}
\P\Big[\|z\|_2^2 \leq (1+4\sqrt{t'})n \Big] \geq 1-\exp(-t' n).
\end{equation}
If we define the event 
\begin{equation}\label{e1}
\mathcal{E}_1:=\Bigg\{\frac{1}{\ts\frac{1}{\|z\|_2^2}\, \Pi_p(z)\ttop M(X) \Pi_p(z) } \leq   \ts\frac{1+4\sqrt{t'}} {(1-2\sqrt{t'})}\ts\frac{n}{\tr(M(X))}\Bigg\}
\end{equation}
then the previous two inequalities imply
\begin{equation}
\begin{split}
\P\big[ \mathcal{E}_1 \big\bracevert M(X) \big] &\geq 1-\exp(-t'\tr(M(X)))-\exp(-t'\cdot n)\\
& \geq 1-2\exp(-t'\tr(M(X))),
\end{split}
\end{equation}
since $\tr(M(X))\leq n$ almost surely.
Next, let $\kappa_5,\kappa_6>0$ be as in the previous lemma, and define the event
\begin{equation}\label{e2}
\mathcal{E}_2:=\Big\{ \ts\frac{n}{\tr(M(X))} \leq \frac{1}{\kappa_5} n^{1-\frac{\gamma}{\eta}}\Big\},
\end{equation}
which has probability $\P(\mathcal{E}_2)\geq 1-\exp(-\kappa_6 n^{1+\frac{\gamma(1-\eta)}{\eta}})$.

We now put these items together. Starting with line~\eqref{quadform}, if we work on the intersection of $\mathcal{E}_1$ and $\mathcal{E}_2$, then for any fixed $i=1,\dots,n$ we have
\small
\begin{equation}
\begin{split}
\P\Big[\ts\frac{1}{v_{\rho}(X;X_i)} \leq \ts\frac{1+4\sqrt{t'}} {(1-2\sqrt{t'})}\ts\frac{1}{\kappa_5} n^{1-\frac{\gamma}{\eta}}\Big]
&\geq \E\Big[ 1_{\mathcal{E}_1} \cdot 1_{\mathcal{E}_2}\Big]\\[0.2cm]
&\geq 1-\E\Big[1_{\mathcal{E}_1^c}+1_{\mathcal{E}_2^c}\Big]\\[0.2cm]
&=1-\E \big[\E\big[1_{\mathcal{E}_1^c}\big\bracevert M(X)\big]\big]-\P(\mathcal{E}_2^c)\\[0.2cm]
&\geq 1-2\E \Big[ \exp(-t'\tr(M(X))\Big]-\P(\mathcal{E}_2^c)\\[0.2cm]
&= 1-2\E \Big[ \exp(-t'\tr(M(X)))\cdot (1_{\mathcal{E}_2}+1_{\mathcal{E}_2^c})\Big]-\P(\mathcal{E}_2^c)\\[0.2cm]
&\geq 1-\exp(-t' \kappa_5 n^{\frac{\gamma}{\eta}})-3\P(\mathcal{E}_2^c)\\[0.2cm]
&\geq 1-\exp(-t'\kappa_5 n^{\frac{\gamma}{\eta}})-3\exp(-\kappa_6 n^{1+\frac{\gamma(1-\eta)}{\eta}})\\[0.2cm]
&\geq  1-4\exp(-\min\{t'\kappa_5,\kappa_6\} \cdot n^{\frac{\gamma}{\eta}})
\end{split}
\end{equation}
\normalsize
where we have used the previous lemma to bound $\P(\mathcal{E}_2^c)$, and also the assumption $\gamma\in (0,1)$ to conclude that  $\frac{\gamma}{\eta} \leq 1+\frac{\gamma(1-\eta)}{\eta}$.
Taking a union bound over $i=1,\dots,n$, proves the claim.\qed
~\\

The last component of Lemma~\ref{propGaussian} is to prove the MSPE inequalities.

\subsection{Proof of the MSPE inequalities}
 
The proof of Theorem~\eqref{ridgemspe} shows that for any realization of $X$ we have
\begin{equation}
\text{mspe}(\hat{\beta}_{\varrho}|X):=\ts\frac{1}{n}\E\Big[\|X(\hatbeta_{\varrho}-\beta)\|_2^2\big\bracevert X \Big] \lesssim n^{-\theta}\|\beta\|_2^2+\ts\frac{1}{n}\displaystyle \sum_{i=1}^{p\wedge n} \big(\ts\frac{l_i}{l_i+\varrho}\big)^2,
\end{equation}
where $l_i=\lambda_i(\frac{1}{n}X\ttop X)$.
Now observe that the second term on the right side matches the expression for $\tr(M(X))$ given in line~\eqref{tracecalc} by replacing $\rho$ with $\varrho$ and multiplying by a factor of $\frac{1}{n}$. Therefore, using Lemma~\ref{traceconc} and recalling $\varrho/n = n^{-\theta}$ shows that there are absolute constants $\kappa_1,\kappa_2,\kappa_3,\kappa_4>0$  such that for all large $n$,
\small
\begin{equation}\label{firstmspe}
\P\Big[ \text{mspe}(\hat{\beta}_{\varrho}|X)\geq \kappa_1\big(n^{-\theta}+ n^{2(\theta-\eta)}\big) \Big] \leq \exp(-\kappa_2 n^{2(1-\eta)+\theta}), \text{ \ \ if \ } \eta \in (0,\ts\frac{1}{2}).
\end{equation}
and
\begin{equation}\label{secondmspe}
\P\Big[ \text{mspe}(\hat{\beta}_{\varrho}|X)\geq \kappa_3\big(n^{-\theta}+ n^{\frac{\theta}{\eta}-1}\big) \Big] \leq \exp(-\kappa_4 n^{1+\frac{\theta(1-\eta)}{\eta}}) , \text{ \ \ \ \ if \ } \eta >\ts\frac{1}{2}.
\end{equation}
\normalsize
In line~\eqref{firstmspe}, the bound $\mspe(\hat{\beta}_{\varrho}|X)$ is optimized when $n^{-\theta}\asymp n^{2(\theta-\eta)}$, which explains the choice $\theta=\frac{2\eta}{3}$. Similarly, in line~\eqref{secondmspe}, the bound is optimized when $n^{-\theta}\asymp n^{\frac{\theta}{\eta}-1}$, which explains the choice $\theta=\frac{\eta}{\eta+1}$. Substituting in these values $\theta$ yields the stated result.\qed

\section{Background results}\label{appmisc}

\subsection{Results on matrices and convexity}\label{convex}
\begin{lemma}\label{frobineq}
Let $M_1\in \R^{k_1\times k_1}$ and $M_2\in \R^{k_1\times k_2}$. Then,
\begin{equation}
\|M_1M_2\|_F\leq \|M_1\|_{\textup{op}}\|M_2\|_F.
\end{equation}
\end{lemma}

\proof Observe that
\begin{equation}
\begin{split}
\|M_1M_2\|_F^2 &= \tr(M_2\ttop M_1\ttop M_1 M_2)\\
&=\tr((M_1\ttop M_1)(M_2M_2\ttop))\\
&\leq \sum_{i=1}^{k_1} \lambda_i(M_1\ttop M_1)\cdot \lambda_i(M_2M_2\ttop)\\
&\leq \|M_1\|\op^2 \sum_{i=1}^{k_1}\lambda_i(M_2M_2\ttop)\\
&=\|M_1\|\op^2 \|M_2\|_F^2,
\end{split}
\end{equation}
where we have used von Neumann's trace inequality (also known as Fan's inequality)~\cite[p.10]{BorweinLewis} in the third line.\qed

\paragraph{A result on convex trace functionals.}

In the following lemma, an interval of the real line refers to any set of the form $(a,b)$,$(a,b]$,$[a,b)$, or $[a,b]$, where $-\infty\leq a\leq b\leq \infty$. We also define $\text{spec}(M)$ to be the set of eigenvalues of a square matrix $M$. The collection of symmetric matrices in $\R^{p \times p}$ is denoted by $\mathbb{S}^{p\times p}$. For a univariate function $\varphi$, the symbol $\tr(\varphi(M))$ denotes $\sum_{i} \varphi(\lambda_i(M))$.

\begin{lemma}\label{convexfunctional}
Let $\mathcal{I}\subset \R$ be an interval, and let $\mathcal{M}\subset \mathbb{S}^{p\times p}$ be a convex set such that $\text{spec}(M)\subset \mathcal{I}$ for all $M\in\mathcal{M}$. Let $\varphi :\mathcal{I}\to \R$ be a convex function. Then, the functional
\begin{equation}
M\mapsto \tr(\varphi(M))
\end{equation}
is convex on $\mathcal{M}$.
\end{lemma}
\noindent A proof may be found in the paper~\cite[Proposition 2]{petz}.

 \paragraph{Operator Jensen inequality.} A function \mbox{$f:\mathbb{S}^{p\times p}\to\mathbb{S}^{p\times p}$} is said to be operator convex if for all $\lambda\in[0,1]$, and all $A,B\in \mathbb{S}^{p\times p}$,
\begin{equation}
f(\lambda A +(1-\lambda)B) \preceq \lambda f(A)+(1-\lambda)f(B),
\end{equation}
where $A\preceq B$ means that $B-A$ is positive semidefinite.
\begin{lemma}[Operator Jensen inequality]\label{opjensen}
Suppose $f:\mathbb{S}^{p\times p} \to \mathbb{S}^{p\times p}$ is operator convex, and let $A$ be a random $\mathbb{S}^{p\times p}$-valued matrix that is integrable. Then,
\begin{equation}
f(\E[A]) \preceq \E[f(A)].
\end{equation}
\end{lemma}
\proof It is enough to show that for all $x\in\R^p$, 
\begin{equation}\label{previousline}
x\ttop f(\E[A])x \leq x\ttop \E[f(A)]x.
\end{equation}
For any fixed $x$, consider the function $g:\mathbb{S}^{p\times p}\to \R$ defined by $g(A)=x\ttop f(A) x$. It is clear that $g$ is a convex function in the usual sense, and so the ordinary version of Jensen's inequality implies $g(\E[A])\leq \E[g(A)]$, which is the same as~\eqref{previousline}.\qed

\subsection{Results on Gaussian vectors and matrices}\label{sec:gaussian}
The following lemma is standard and is often referred to as the Gaussian concentration inequality~\cite{boucheron}.
\begin{lemma}\label{GaussianConc}
Let $Z\in \R^p$ be a standard Gaussian vector and let $f:\R^p\to \R$ be an $L$-Lipschitz function with respect to the $\ell_2$ norm. Then for all $t>0$,
\begin{equation}
\P\Big(|f(Z)-\E[f(Z)]|\geq t\Big)\leq 2\exp\big(\ts\frac{-t^2}{2L^2}\big).
\end{equation}¥
\end{lemma}

Next, we give a formula for the expected square of a Wishart matrix.
\begin{lemma}\label{wishartsquared}
Let $X\in \R^{n\times p}$ have rows drawn i.i.d. from $N(0,\Sigma)$, and let $\hat{\Sigma} = \frac{1}{n}X\ttop X$. Then,
$$\E[ \hat{\Sigma}^2 ] = (1+\ts\frac{1}{n}) \Sigma^2+\ts\frac{\tr(\Sigma)}{n} \Sigma.$$
\end{lemma}

\proof Write $\hat{\Sigma} = \ts\frac{1}{n}\sum_{i=1}^n X_iX_i\ttop$ where $X_i\ttop \in \R^p$ is the $i$th row of $X$. If we put $M_i = X_iX_i\ttop$, then
$$\hat{\Sigma}^2 = \ts\frac{1}{n^2}\sum_{i=1}^n M_i^2+ \ts\frac{1}{n^2}\sum_{i\neq j}M_i M_j. $$
Clearly, the $M_i$ are independent with $\E[M_i] = \Sigma$ for all $i$, and so
$$\E[\hat{\Sigma}^2] = \ts\frac{1}{n}\E[M_1^2]+ \ts\frac{2}{n^2}\binom{n}{2}\Sigma^2.$$
It remains to compute $\E[M_1^2]$. Write $X_i= \Sigma^{1/2}Z_i$ where $Z_i$ is a standard Gaussian vector in $\R^p$, and let $\Sigma^{1/2} = U \Lambda^{1/2} U\ttop$ be a spectral decomposition of $\Sigma^{1/2}$ where $U\in \R^{p\times p}$ is orthogonal, and $\Lambda$ is diagonal with $\Lambda_{ii}=\lambda_i(\Sigma)$. By the orthogonal invariance of the normal distribution, $X_i \de U \Lambda^{1/2}Z_i$, and so
\begin{equation}\label{interm}
\E[M_1^2] = \E[X_1X_1\ttop X_1X_1\ttop] = U\Lambda^{1/2} \E\Big[ Z_1Z_1\ttop \Lambda Z_1 Z_1\ttop \Big]\Lambda^{1/2} U\ttop.
\end{equation}
Define the matrix $\tilde{M}_1:=Z_1Z_1\ttop \Lambda Z_1Z_1\ttop$. It is straightforward to verify that $\E[\tilde{M}_1]$ is diagonal, and its $j$th diagonal entry is
$$\E[\tilde{M}_1]_{jj}
=\tr(\Sigma)+2\lambda_j.$$
Therefore,
$$\E[\tilde{M}_1] = \tr(\Sigma)I_{p\times p} +2\Lambda,$$
and combining this with line~\eqref{interm} gives
$$\E[M_1^2] = \tr(\Sigma) \Sigma+ 2\Sigma^2.$$\qed

\paragraph{Signed s.v.d.}The following lemma describes the factors of an s.v.d. of a standard Gaussian matrix $Z\in\R^{n\times p}$ with $n\geq p$. 
 To make the statement of the lemma more concise, we define the term \emph{signed s.v.d.} below. This is merely a particular form of the s.v.d. that ensures uniqueness. Specifically, if $Z\in\R^{n\times p}$ is a full rank matrix with $n\geq p$, then the signed s.v.d. of $Z$ is given by
\begin{equation}
Z=HLG\ttop,
\end{equation}
where $H\in \R^{n\times p}$ has orthonormal columns, the matrix $L\in\R^{p\times p}$ is diagonal with $L_{11}\geq L_{22}\geq\cdots \geq L_{pp}>0$, and $G\in\R^{p\times p}$ is orthogonal with its first row non-negative, i.e. $G_{1i}\geq 0$ for all $i=1,\dots,p$. It is a basic fact from linear algebra that the signed s.v.d. of any full rank matrix in $\R^{n\times p}$ exists and is unique~\cite[Lemma 7.3.1]{hornjohnson}. 
This fact applies to Gaussian matrices, since they are full rank with probability 1.

\begin{lemma}\label{svdlemma}
Suppose $n\geq p$, and let $Z\in\R^{n\times p}$ be a random matrix with entries drawn i.i.d. from $N(0,1)$. Let
\begin{equation}
Z= HLG\ttop,
\end{equation}
be the unique signed s.v.d.~of $Z$ as defined above.
Then, the matrices $H$, $L$ and $G$ are independent. Furthermore, if $H_i\ttop\in\R^p$ denotes the $i$th row of $H$, then for each $i=1,\dots,n$, the marginal law of $H_i\ttop$ is given by
\begin{equation}\label{hrep}
H_i\ttop \overset{\mathcal{L}}{=} \ts\frac{1}{\|z\|_2}\cdot\Pi_p(z),
\end{equation}
where $z\in\R^n$ is a standard Gaussian vector, and $\Pi_p(z)$ is the projection operator onto the first $p$ coordinates, i.e. $\Pi_p(z)= (z_1,\dots,z_p)$.
\end{lemma}
\proof We first argue that $H$, $L$, and $G$ are independent, and then derive the representation for $H_i\ttop$ in the latter portion of the proof. 

Due to the fact that the transformation $Z\mapsto (H,L,G)$ is invertible, it is possible to obtain the joint density of $(H,L,G)$ from the density of $Z$ by computing the matrix Jacobian of the factorization $Z=HLG$. (See the references~\cite{muirhead},~\cite{mathai}, and~\cite{edelmanRMT} for more background on Jacobians of matrix factorizations.) To speak in more detail about the joint density, let $\mathbb{V}^{n\times p}$ denote the Stiefel manifold of $n\times p$ matrices with orthonormal columns. Also, let $\mathbb{D}^{p\times p}$ denote the set of $p\times p$ diagonal matrices, and $\mathbb{O}^{p\times p}$ the set of orthogonal $p\times p$ matrices. The subset of $\mathbb{O}^{p\times p}$ with non-negative entries in the first row will be denoted by $\mathbb{O}^{n\times p}_+$.

Let \mbox{$f_{H,L,G}:\mathbb{V}^{n\times p} \times  \mathbb{D}^{p\times p}\times \mathbb{O}^{p\times p}_+\to [0,\infty)$} denote the joint density of $(H,L,G)$, where the base measure is the product of Haar measure on $\mathbb{V}^{n\times p}$, Lebesgue measure on $\mathbb{D}^{p\times p}$, and Haar measure on $\mathbb{O}^{p\times p}$ restricted to $\mathbb{O}^{p\times p}_+$. (See the book~\cite{chikuse} for background on these measures). From lines 8.8-8.10 in the paper~\cite{james1954normal}, it is known that such a density exists, and that $f_{H,L,G}$ factors according to
\begin{equation}\label{factordensity}
f_{H,L,G}(h,l,g) = f_H(h)\cdot f_L(l)\cdot f_G(g),
\end{equation}
where the function $f_H$ denotes the density of $H$ with respect to Haar measure on $\mathbb{V}^{n\times p}$, the function $f_L$ denotes the density of $L$ with respect to Lebesgue measure on $\mathbb{D}^{p\times p}$, and the function $f_G$ denotes the density of $G$ with respect to the restriction of Haar measure to $\mathbb{O}^{p\times p}_+$. 
This proves that $H$, $L$, and $G$ are independent.

We now prove the representation~\eqref{hrep}. From line 8.10 in the paper~\cite{james1954normal}, it is known that $f_H$ is constant on $\mathbb{V}^{n\times p}$. In other words, the matrix $H$ follows the Haar distribution on $\mathbb{V}^{n\times p}$.
Consequently, Theorem 2.2.1(ii) in the book \cite{chikuse} implies that the rows $H_i\ttop$ can be represented as
\begin{equation}
H_i\ttop \overset{\mathcal{L}}{=} \Pi_p(J_i\ttop)
\end{equation}
where $J_i\ttop\in \R^n$ is the $i$th row of a Haar-distributed random matrix in $\mathbb{O}^{n\times n}$. Furthermore, the rows $J_i\ttop$ are uniformly distributed on the unit sphere on $\R^n$, and hence can be represented as $z/\|z\|_2$, where $z\in\R^n$ is a standard Gaussian vector.\qed

\bibliographystyle{unsrt}
{\bibliography{BootstrapBib_3july2016_arxiv.bib}}

\begin{thebibliography}{10}

\bibitem{zhang2014confidence}
C.-H. Zhang and S.~S. Zhang.
\newblock Confidence intervals for low dimensional parameters in high
  dimensional linear models.
\newblock {\em Journal of the Royal Statistical Society: Series B},
  76(1):217--242, 2014.

\bibitem{javanmard2013hypothesis}
A.~Javanmard and A.~Montanari.
\newblock Hypothesis testing in high-dimensional regression under the
  {G}aussian random design model: Asymptotic theory.
\newblock {\em arXiv preprint arXiv:1301.4240}, 2013.

\bibitem{javanmard2013confidence}
A.~Javanmard and A.~Montanari.
\newblock Confidence intervals and hypothesis testing for high-dimensional
  regression.
\newblock {\em arXiv preprint arXiv:1306.3171}, 2013.

\bibitem{buhlmannridge}
P.~B{\"u}hlmann.
\newblock Statistical significance in high-dimensional linear models.
\newblock {\em Bernoulli}, 19(4):1212--1242, 2013.

\bibitem{vandegeerci2013}
S.~van~de Geer, P.~B{\"u}hlmann, and Y.~Ritov.
\newblock On asymptotically optimal confidence regions and tests for
  high-dimensional models.
\newblock {\em arXiv preprint arXiv:1303.0518}, 2013.

\bibitem{lee2013}
J.~D. Lee, D.~L. Sun, Y.~Sun, and J.~E. Taylor.
\newblock Exact inference after model selection via the lasso.
\newblock {\em arXiv preprint arXiv:1311.6238}, 2013.

\bibitem{lahiribootstrap2013}
A.~Chatterjee and S.~N. Lahiri.
\newblock Rates of convergence of the adaptive lasso estimators to the oracle
  distribution and higher order refinements by the bootstrap.
\newblock {\em The Annals of Statistics}, 41(3):1232--1259, 2013.

\bibitem{binyuridge2013}
H.~Liu and B.~Yu.
\newblock Asymptotic properties of lasso+mls and lasso+ridge in sparse
  high-dimensional linear regression.
\newblock {\em Electronic Journal of Statistics}, 7:3124--3169, 2013.

\bibitem{chernozhukovmultiplier}
V.~Chernozhukov, D.~Chetverikov, and K.~Kato.
\newblock Gaussian approximations and multiplier bootstrap for maxima of sums
  of high-dimensional random vectors.
\newblock {\em The Annals of Statistics}, 41(6):2786--2819, 2013.

\bibitem{lehmann}
E.~L. Lehmann and J.~P. Romano.
\newblock {\em Testing Statistical Hypotheses}.
\newblock Springer, 2005.

\bibitem{freedman1981}
D.~A. Freedman.
\newblock Bootstrapping regression models.
\newblock {\em The Annals of Statistics}, 9(6):1218--1228, 1981.

\bibitem{bickel1983}
P.~J. Bickel and D.~A. Freedman.
\newblock Bootstrapping regression models with many parameters.
\newblock In {\em Festschrift for Erich L. Lehmann}, pages 28--48. Wadsworth,
  1983.

\bibitem{draper1998}
N.~R. Draper and H.~Smith.
\newblock {\em Applied Regression Analysis}.
\newblock Wiley-Interscience, 1998.

\bibitem{bickel1981}
P.~J. Bickel and D.~A. Freedman.
\newblock Some asymptotic theory for the bootstrap.
\newblock {\em The Annals of Statistics}, pages 1196--1217, 1981.

\bibitem{bobkovledoux}
S.~Bobkov and M.~Ledoux.
\newblock {\em One-dimensional empirical measures, order statistics, and
  Kantorovich transport distances}.
\newblock preprint, 2014.

\bibitem{tsybakov2009}
A.~B. Tsybakov.
\newblock {\em Introduction to Nonparametric Estimation}.
\newblock Springer, 2009.

\bibitem{wassermanallnp}
L.~Wasserman.
\newblock {\em All of Nonparametric Statistics}.
\newblock Springer, 2006.

\bibitem{LaurentMassart}
B.~Laurent and P.~Massart.
\newblock Adaptive estimation of a quadratic functional by model selection.
\newblock {\em Annals of Statistics}, 28(5):1302--1338, 2000.

\bibitem{horntopics}
R.~A. Horn and C.~R. Johnson.
\newblock Topics in matrix analysis.
\newblock {\em Cambridge University Press}, 1991.

\bibitem{bhatia}
R.~Bhatia.
\newblock {\em Matrix Analysis}.
\newblock Springer, 1997.

\bibitem{hornjohnson}
R.~A. Horn and C.~A. Johnson.
\newblock {\em Matrix Analysis}.
\newblock Cambridge University Press, 22nd printing edition, 2009.

\bibitem{muirhead}
R.~J. Muirhead.
\newblock {\em Aspects of Multivariate Statistical Theory}.
\newblock John Wiley \& Sons, 1982.

\bibitem{BorweinLewis}
J.~M. Borwein and A.~S. Lewis.
\newblock {\em Convex Analysis and Nonlinear Optimization Theory and Examples}.
\newblock CMS Bookks in Mathematics. Canadian Mathematical Society, 2000.

\bibitem{petz}
D.~Petz.
\newblock A survey of certain trace inequalities.
\newblock {\em Functional analysis and operator theory}, 30:287--298, 1994.

\bibitem{boucheron}
S.~Boucheron, G.~Lugosi, and P.~Massart.
\newblock {\em Concentration Inequalities: A Nonasymptotic Theory of
  Independence}.
\newblock Oxford University Press, 2013.

\bibitem{mathai}
A.~M. Mathai.
\newblock {\em Jacobians of Matrix Transformations and Functions of Matrix
  Argument}.
\newblock World Scientific, 1997.

\bibitem{edelmanRMT}
A.~Edelman and N.~R. Rao.
\newblock Random matrix theory.
\newblock {\em Acta Numerica}, 14:233--297, 2005.

\bibitem{chikuse}
Y.~Chikuse.
\newblock {\em Statistics on Special Manifolds}, volume 174.
\newblock Springer Science \& Business Media, 2003.

\bibitem{james1954normal}
A.~T. James.
\newblock Normal multivariate analysis and the orthogonal group.
\newblock {\em The Annals of Mathematical Statistics}, pages 40--75, 1954.

\end{thebibliography}


\end{document}